\documentclass[preprint,12pt]{elsarticle}

\usepackage{amsfonts}
\usepackage{mathtools}
\usepackage{enumitem}
\usepackage{booktabs}
\usepackage{multirow}
\usepackage{mdframed}
\usepackage{subcaption}
\usepackage[section]{placeins}
\usepackage{xcolor}
\usepackage{url}

\newcommand{\R}{\mathbb{R}}
\newcommand{\norm}[1]{\|#1\|}
\newcommand{\abs}[1]{|#1|}
\newcommand{\ip}[2]{\langle #1, #2 \rangle}

\newcommand{\sign}{\operatorname{sign}}

\newcommand{\argmin}{\mathop{\mathrm{arg\,min}}}
\newcommand{\diag}{\operatorname{diag}}

\newcommand{\bd}{\mathbf{d}}
\newcommand{\bp}{\mathbf{p}}

\newcommand{\bb}{\mathbf{b}}
\newcommand{\bx}{\mathbf{x}}
\newcommand{\bu}{\mathbf{u}}
\newcommand{\bz}{\mathbf{z}}

\newcommand{\by}{\mathbf{y}}

\newcommand{\ba}{\mathbf{a}}
\newcommand{\bv}{\mathbf{v}}
\newcommand{\bzero}{\mathbf{0}}
\newcommand{\bone}{\mathbf{1}}

\begin{document}

\begin{frontmatter}

\title{Exact Total Variation Minimizers as Non-Oscillatory Limiters in High-Order Methods for Conservation Laws}

\author[uiuc]{Gabriel P. Langlois}
\ead{gp42@illinois.edu}
\author[brown]{J\'er\^ome Darbon}
\ead{jerome\_darbon@brown.edu}
\author[purdue]{Rongjie Lai}
\ead{lairj@purdue.edu}
\author[brown]{Chi-Wang Shu}
\ead{chi-wang\_shu@brown.edu}
\author[purdue]{Xiangxiong Zhang\corref{cor1}}
\ead{zhan1966@purdue.edu}
\cortext[cor1]{Corresponding author}

\address[uiuc]{Department of Mathematics,
University of Illinois Urbana-Champaign, Urbana, IL 61801, USA}
\address[brown]{Division of Applied Mathematics,
Brown University, Providence, RI 02912, USA}
\address[purdue]{Department of Mathematics,
Purdue University, West Lafayette, IN 47907, USA}

\begin{abstract}
High-order numerical methods for conservation laws can generate spurious oscillations near discontinuities.  We propose to suppress these oscillations by post-processing the numerical solution with a total variation (TV) denoising step that reduces the total variation of the nodal values. For a general analysis operator, the resulting discrete minimization problem need not satisfy the submodularity property that existing exact max-flow algorithms require, so we instead compute the TV minimizer exactly using a differential inclusion algorithm. In one dimension, the differential inclusion algorithm computes the TV minimizer in finitely many steps, requires no tuning of algorithmic parameters, and preserves the total mass. In two dimensions, we apply the 1D algorithm dimension by dimension. We test the method as a post-processing limiter for Fourier pseudospectral methods and a fifth-order finite difference scheme applied to scalar conservation laws, compressible Euler equations, and the two-dimensional incompressible Euler equations.
\end{abstract} 


\end{frontmatter}

\section{Introduction}
\label{sec:intro}

\subsection{Oscillations in high-order methods for conservation laws}
\label{sec:intro-oscillations}

Consider the scalar conservation law
\begin{equation}\label{eq:intro-cl}
  u_t + f(u)_x = 0
\end{equation}
with suitable initial data.
It is well known that solutions of~\eqref{eq:intro-cl} develop
discontinuities in finite time even from smooth initial data, and that
linear high-order methods
for smooth solutions generate spurious oscillations near such
discontinuities.
For spectral methods~\cite{gottlieb1977,hesthaven2007spectral}, the global
Fourier or polynomial basis produces Gibbs oscillations that are $O(1)$
in amplitude and that pollute the entire computational domain.
See~\cite{Shu1998,shu2009} for non-oscillatory high-order finite difference and finite volume
methods such as ENO and WENO schemes. For discontinuous Galerkin (DG)
methods~\cite{cockburn2001,hesthaven2007}, limiters are needed to treat troubled cells.

Existing remedies are usually designed within a specific discretization framework. For spectral methods, common approaches include exponential filters~\cite{vandeven1991}, spectral viscosity~\cite{tadmor1989}, and Gegenbauer reconstruction~\cite{gottliebshu1997}. For finite volume and DG methods, common approaches include TVD limiters~\cite{harten1983}, WENO reconstruction~\cite{shu2009,JiangShu1996}, invariant-domain-preserving (IDP) limiters~\cite{ZhangShu2010,ZhangShu2012fd,guermond2018}, flux-corrected transport (FCT)~\cite{borisbook1973,zalesak1979}, and optimization-based limiters~\cite{guba2014,bochev2012constrained}. We refer to~\cite{wu2025high} for a comprehensive survey of IDP schemes. The bound-preserving sweeping technique~\cite{liu2017sweep} enforces pointwise bounds while conserving a weighted average, but it does not directly reduce the total variation.

\subsection{TV regularization}
\label{sec:intro-tv}
Given a numerical solution $\bb \in \R^n$ contaminated by oscillations,
the TV denoising problem seeks
\begin{equation}\label{eq:intro-tv}
  \bx^* = \argmin_{\bx \in \R^n}\;
  \|D\bx\|_1 + \frac{1}{2\lambda}\|\bx - \bb\|_2^2,
\end{equation}
where $D \in \R^{(n-1)\times n}$ is the forward difference matrix, i.e., $(D \bx)_i=x_{i+1}-x_i$, and
$\lambda > 0$ is a chosen parameter to control the regularization strength. The minimizer of~\eqref{eq:intro-tv} preserves the total mass of the input: $\sum_i x^*_i = \sum_i b_i$. This follows from $\ker(D) = \operatorname{span}\{\bone\}$ (see
Section~\ref{sec:reduction}) and ensures that the TV limiter defined by~\eqref{eq:intro-tv} conserves total mass. TV regularization was introduced for image denoising~\cite{rudin1992,chambolle2004}. In this paper, we use~\eqref{eq:intro-tv} as a limiter to reduce the total variation of the nodal values.

\subsection{Computational cost of TV solvers}
\label{sec:intro-cost}

In a PDE context where the TV filter may be invoked at every time step,
iterative solver cost can dominate the computation. The key focus of this paper is to introduce an efficient exact solver for~\eqref{eq:intro-tv}.

\begin{table}[!htbp]
  \centering
  \caption{Solver comparison for selected test problems.
    Steps and DI time are for the differential inclusion
    algorithm~\cite{langlois2025} (C++ implementation).
    CD (C++ implementation) time is for the Chambolle--Darbon parametric
    max-flow solver~\cite{chambolle2009,hochbaum2001}.
    All times are online solve times on a MacBook Pro with Apple M1 chip
    (offline precomputation of $(D^\top D)^\dagger D^\top$ excluded).
    Times are averages over $R$ runs, where $R$ is chosen so
    that the total elapsed time exceeds one second.}
  \label{tab:di-summary}
  \smallskip\footnotesize
  \begin{tabular}{@{}llrrrrl@{}}
    \hline
    Test problem & $n$ & $\lambda$ & DI Steps & DI (s) & CD (s) & Sec. \\
    \hline
    Gaussian noise, pw.\ const.
      & 200  & 1.0  & 13  & $4.5\text{e-}5$ & $3.6\text{e-}5$ & \S\ref{sec:tv-gaussian} \\
      & 1000 & 1.0  & 59  & $2.5\text{e-}3$ & $3.1\text{e-}4$ & \S\ref{sec:tv-gaussian} \\[3pt]
    Spectral Gibbs, pw.\ const.\
      & 200  & 0.3  & 3   & $1.8\text{e-}5$ & $4.6\text{e-}5$ & \S\ref{sec:tv-spectral} \\
      & 2000 & 0.3  & 3   & $1.5\text{e-}3$ & $3.2\text{e-}3$ & \S\ref{sec:tv-spectral} \\[3pt]
    Spectral Gibbs, pw.\ smooth
      & 200  & 0.08 & 110 & $3.4\text{e-}3$ & $3.5\text{e-}5$ & \S\ref{sec:tv-spectral2} \\
      & 1000 & 0.08 & 562 & $4.2\text{e-}1$ & $2.4\text{e-}4$ & \S\ref{sec:tv-spectral2} \\[3pt]
    FD5 Burgers post-proc.
      & 256  & 0.02 & 244 & $4.1\text{e-}2$ & $4.2\text{e-}5$ & \S\ref{sec:pde-burgers} \\
    \hline
  \end{tabular}
\end{table}

For the one-dimensional problem~\eqref{eq:intro-tv}, there are several specialized exact or fast solvers, including the taut string algorithm~\cite{davies2001}, the direct algorithm of Condat~\cite{condat2013}, and LARS/homotopy path-following for the generalized LASSO~\cite{tibshirani2011,efron2004}. The 1D TV problem reduces to a LASSO problem through the change of variables $\bz = D\tilde\bx$, where $\tilde\bx = \bx - \operatorname{mean}(\bx)\,\bone$ is the mean-subtracted part of $\bx$; see Section~\ref{sec:reduction}. The taut string algorithm is restricted to the 1D TV problem, and the Condat algorithm extends to the 1D fused LASSO (TV plus $\ell^1$ sparsity), with $O(n)$ complexity in practice~\cite{condat2013}.

An earlier exact solver for the TV denoising problem~\eqref{eq:intro-tv}, applicable when
the problem satisfies a certain submodularity property, is the parametric max-flow algorithm~\cite{hochbaum2001,darbon2005fast,darbon2005composants,darbon2006image,darbon2006imageII,chambolle2009}. This algorithm reformulates the problem as a sequence of binary min-cut problems on the same graph, solved simultaneously by a single parametric max-flow computation. For the 1D TV problem, this reduces to a path graph.

The recent differential inclusion (DI) algorithm for the LASSO problem by Langlois and Darbon~\cite{langlois2025} exploits the polyhedral structure of the dual problem by tracing the piecewise continuous trajectory of its associated differential inclusions exactly, terminating in finitely many steps. Their algorithm is a novel extension of the differential inclusions algorithm proposed in~\cite{tendero2021algorithm} for solving the equality-constrained basis pursuit problem efficiently and exactly, numerically up to machine precision.

Table~\ref{tab:di-summary} compares the two solvers: for test cases with few breakpoints, an optimized C++ implementation of the DI algorithm is~\emph{comparable to or faster} than the parametric max-flow algorithm (e.g., for the 3-step piecewise constant cases the DI solver is roughly $2$--$3\times$ faster at both $n = 200$ and $n = 2000$), while for cases with many
breakpoints (tens of steps and up) the max-flow algorithm is substantially faster. The max-flow algorithm requires a certain submodularity property, which does not hold for other analysis operators $D$ such as those in Table~\ref{tab:fd4-summary} or~\ref{app:operators}.  The DI solver solves the LASSO for any measurement matrix, so it applies to any surjective analysis operator~$D$, e.g., see Table~\ref{tab:fd4-summary} and examples (db1 and db4 wavelets) in~\ref{app:operators}.

\begin{table}[!htbp]
  \centering
  \caption{CPU time comparison of MATLAB and C++ implementations of the
    DI algorithm with
    $D$ replaced by a \emph{fourth-order accurate} finite difference
    approximation to $du/dx$ (interior rows: $[1,-8,0,8,-1]/12$,
    boundary rows: $[-1,1]$).
    All times are online solve times on a MacBook Pro with Apple M1 chip.\label{tab:fd4-summary}}
  \smallskip\footnotesize
  \begin{tabular}{llrrrr}
    \hline
    Test problem & $n$ & $\lambda$ & DI Steps & DI MATLAB (s) & DI C++ (s) \\
    \hline
    Gaussian noise, pw.\ const.
      & 200  & 1.0  &  32 & $1.0\text{e-}2$ & $3.0\text{e-}4$ \\
      & 1000 & 1.0  & 138 & $2.8\text{e-}1$ & $2.7\text{e-}2$ \\[3pt]
    Spectral Gibbs, pw.\ const.\
      & 200  & 0.5  &  24 & $2.0\text{e-}3$ & $2.0\text{e-}4$ \\
      & 2000 & 0.5  & 282 & $3.1\text{e+}0$ & $3.8\text{e-}1$ \\[3pt]
    Spectral Gibbs, pw.\ smooth
      & 200  & 0.5  & 175 & $5.4\text{e-}2$ & $3.3\text{e-}3$ \\
      & 1000 & 0.5  & 1441 & $2.4\text{e+}1$ & $1.3\text{e+}0$ \\
    \hline
  \end{tabular}
\end{table}

For a general operator~$D$, popular iterative methods for~\eqref{eq:intro-tv} include standard iterative solvers such as the Chambolle--Pock primal-dual hybrid gradient (PDHG) method~\cite{chambolle2011}, and the alternating direction method of multipliers (ADMM), which can be interpreted through Douglas--Rachford splitting~\cite{lionsmercier1979,demanet2016,zhang2025}. Their performance depends on algorithmic parameters that require problem-dependent tuning.

For the simpler scalar bound-preserving optimization limiter with a conservation constraint, Bradley et al.~\cite{bradley2019communication} showed that the $\ell^2$ minimizer is also an $\ell^1$ minimizer and that an $\ell^1$ minimizer can be constructed explicitly by clipping to bounds and redistributing the mass deficit, though $\ell^1$ minimizers are generally non-unique~\cite{liu2024opt}. Such a  construction does not extend to the $\ell^1$-penalized LASSO arising from TV denoising.

\subsection{Main results}
\label{sec:intro-results}

We apply the DI
algorithm~\cite{langlois2025} developed by Langlois and Darbon to the TV denoising
problem~\eqref{eq:intro-tv} and use it as a non-oscillatory limiter
for numerical conservation law solvers. The algorithm computes the minimizer of~\eqref{eq:intro-tv} in finitely
many steps, with the step count depending on signal complexity
rather than grid size~$n$ for
clean piecewise-constant data, such as the truncated Fourier sums of
step functions in Section~\ref{sec:tv-spectral}; for noisy
piecewise-constant data (Section~\ref{sec:tv-gaussian}) and for
piecewise-smooth signals the step count grows with~$n$
(Section~\ref{sec:tv-spectral2}).
The matrix $(D^\top D)^\dagger D^\top$ (where $\dagger$ denotes the
Moore--Penrose pseudoinverse) is formed once as an offline
precomputation in $O(n^2\log n)$ work using the known eigendecomposition of $D^\top D$
(the discrete Neumann Laplacian); for the general fourth-order
operator of Table~\ref{tab:fd4-summary}, $M = D^\dagger$ is instead
formed by a dense pseudoinverse offline. The times reported below include
only the online solve. Table~\ref{tab:di-summary} summarizes the
solver performance for selected test problems.

Table~\ref{tab:fd4-summary} demonstrates the generality of the
DI algorithm when we replace the first-difference operator~$D$
with a fourth-order accurate finite difference approximation to $du/dx$
(5-point interior stencil). The max-flow algorithm cannot handle this operator because its graph-cut formulation requires each row of~$D$ to couple at most two variables with the sign pattern that makes the pairwise terms submodular.
The DI algorithm applies without modification: only the
offline pseudoinverse $(D^\top D)^\dagger D^\top$ changes.
In Table~\ref{tab:fd4-summary}, we also list the comparison of MATLAB and C++ implementations of the DI method.

As a limiter, the TV denoising step~\eqref{eq:intro-tv} is applied to
the nodal values of the numerical solution at each time step or as a
post-processing filter.
We test the limiter on scalar convection equations using Fourier
pseudospectral discretizations, and on the 
compressible Euler equations using a fifth-order finite difference
scheme.
In two dimensions, we apply the 1D TV limiter dimension-by-dimension.

\subsection{Organization}
\label{sec:intro-org}
The rest of the paper is organized as follows. 
Section~\ref{sec:formulation} formulates TV denoising as a
non-oscillatory limiter, derives the reduction to LASSO, describes the
dimension-by-dimension extension to 2D, and reviews PDHG and ADMM.
Section~\ref{sec:algorithm} presents the DI algorithm and tests it on a few examples as a demonstration of its performance.
Section~\ref{sec:pde} applies the TV limiter to Fourier pseudospectral
and fifth-order finite difference computations of conservation laws.
The concluding remarks are given in Section~\ref{sec:conclusions}.

\section{TV Limiter Formulation}
\label{sec:formulation}

\subsection{TV denoising as a non-oscillatory limiter}

\begin{figure}[!htbp]
  \centering
  \begin{subfigure}[t]{0.48\textwidth}
    \centering
    \includegraphics[width=\textwidth]{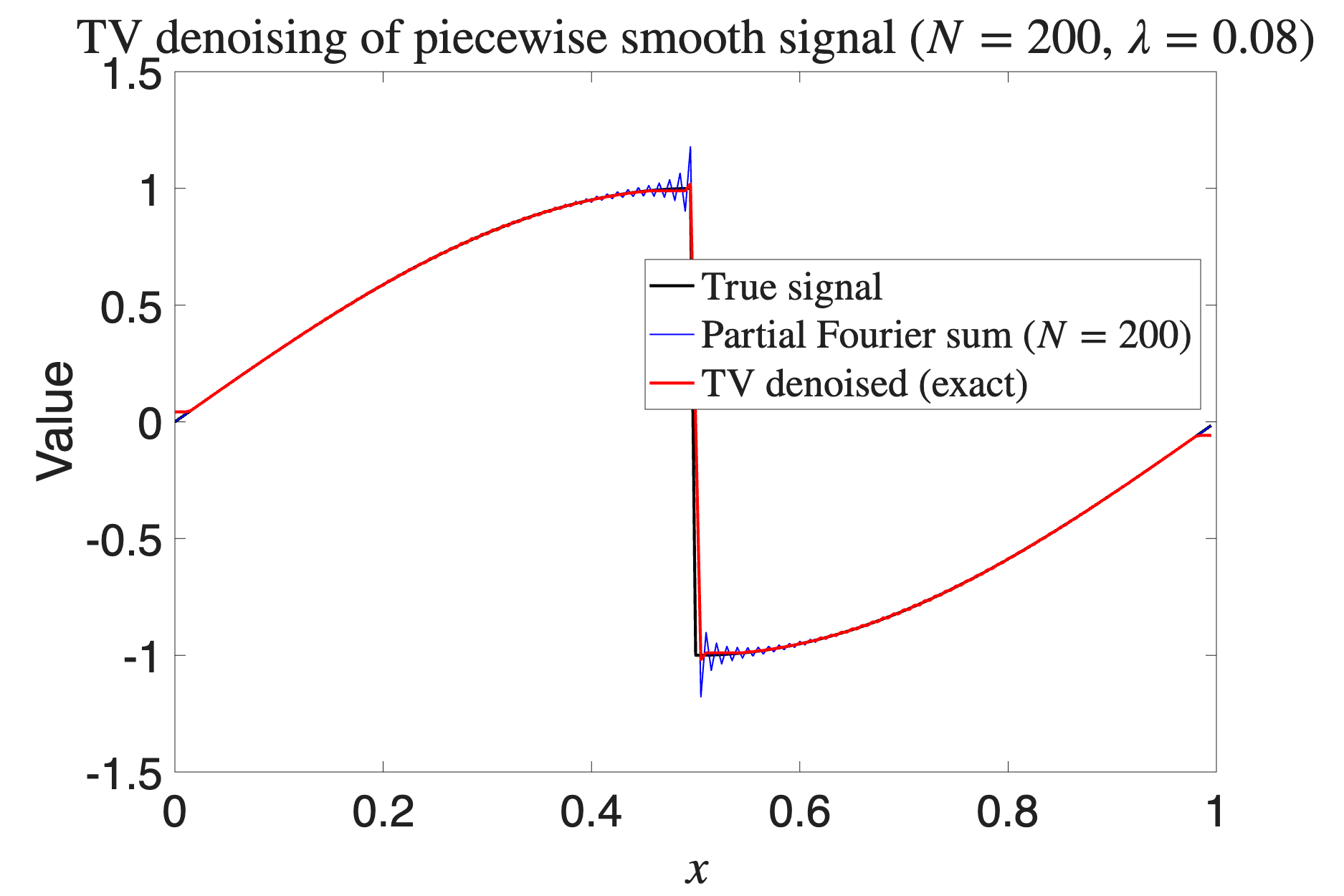}
    \caption{Piecewise smooth signal ($N=200$ Fourier modes,
    $\lambda=0.08$). Blue: partial Fourier sum. Black: true piecewise
    sinusoidal signal. Red: TV denoised.}
    \label{fig:tv-spectral2-signal}
  \end{subfigure}
  \hfill
  \begin{subfigure}[t]{0.48\textwidth}
    \centering
    \includegraphics[width=\textwidth]{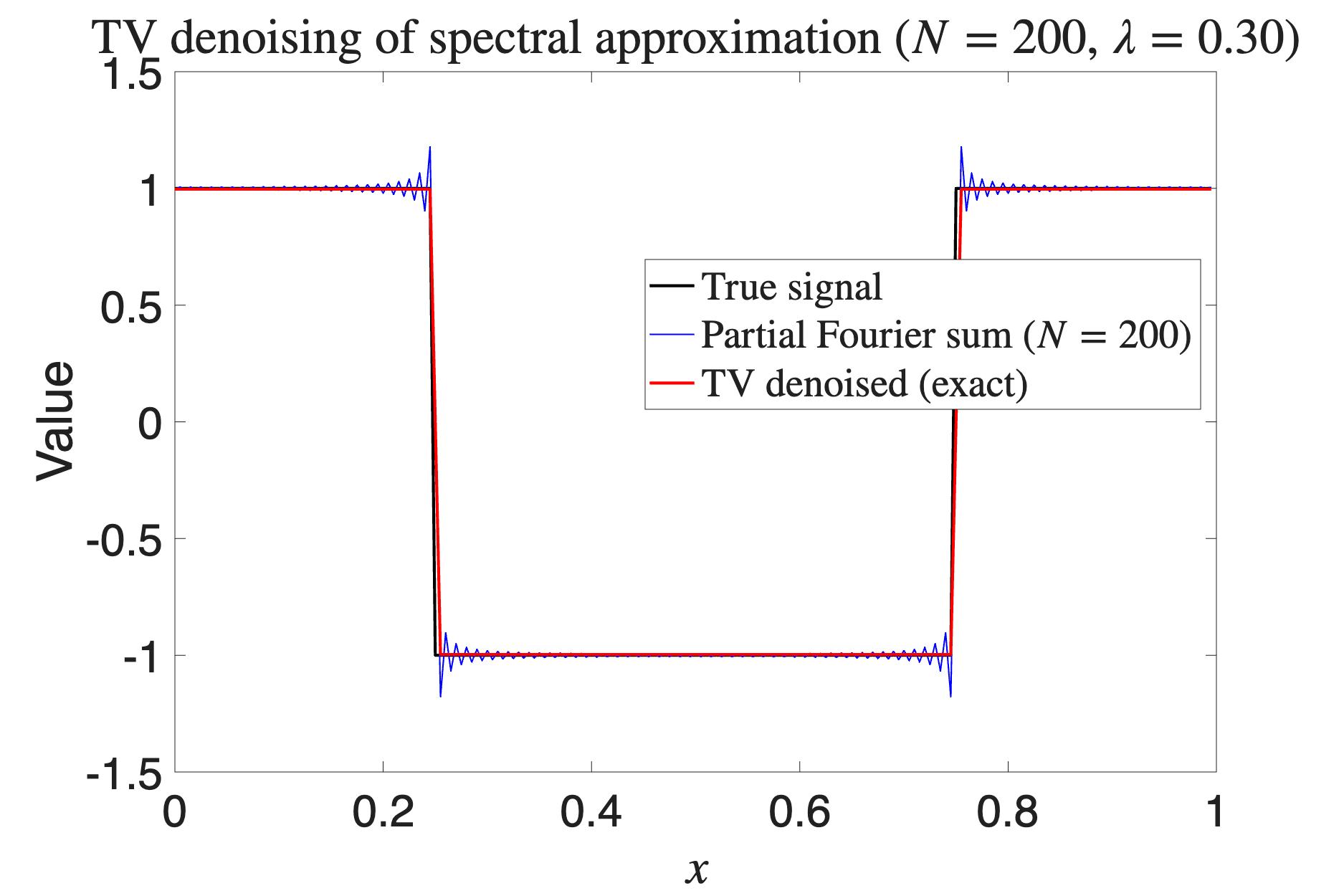}
    \caption{Piecewise constant signal ($N=200$ Fourier modes,
    $\lambda=0.3$). Blue: partial Fourier sum. Black: true step function.
    Red: TV denoised.}
    \label{fig:tv-spectral-signal}
  \end{subfigure}
  \caption{TV denoising as a non-oscillatory filter or limiter for removing
  Gibbs oscillations from truncated Fourier approximations.
  Left: piecewise smooth signal with staircasing artifact visible
  but mild.  Right: piecewise constant signal recovered to plotting
accuracy.}
  \label{fig:tv-denoising-demo}
\end{figure}

We post-process the numerical solution by solving the TV denoising problem~\eqref{eq:intro-tv}, where $D$ is a
first-order finite difference operator and $\lambda > 0$ controls the
strength of regularization. We refer to this approach as a \emph{TV limiter}.
Figure~\ref{fig:tv-denoising-demo} shows TV denoising applied to two
spectral approximations with Gibbs oscillations.

The  forward difference matrix $D \in \R^{(n-1) \times n}$ is
\begin{equation*}
  D =
  \begin{pmatrix}
    -1 & 1 & & & \\
       & -1 & 1 & & \\
       & & \ddots & \ddots & \\
       & & & -1 & 1
  \end{pmatrix},
  \quad
  (D\bx)_i = x_{i+1} - x_i, \quad i = 1,\ldots,n-1,
\end{equation*}
with $\ker(D) = \operatorname{span}\{\bone\}$.
The $\ell_1$ penalty promotes piecewise constant structure,
suppressing oscillations while preserving sharp jumps.
Since $\ker(D) = \operatorname{span}\{\bone\}$, the TV
penalty is blind to the mean of~$\bx$, so the quadratic fidelity
term forces $\sum_i x^*_i = \sum_i b_i$, preserving the total mass
exactly.   

\subsection{Reduction to LASSO}
\label{sec:reduction}

For any vector $\bv \in \R^n$, write its mean-subtracted part as
$\tilde\bv = \bv - \bar{v}\bone$ with $\bar v = \operatorname{mean}(\bv)$,
so that $\tilde\bv \perp \bone$.
Since $D\bone = \bzero$, we decompose $\bx = \bar{x}\bone + \tilde\bx$,
so that the TV term depends only on $\tilde\bx$. The full objective
then splits into two independent pieces:
\[
  \norm{D\bx}_1 + \frac{1}{2\lambda}\norm{\bx - \bb}_{2}^{2}
  = \underbrace{\frac{n}{2\lambda}(\bar{x} - \bar{b})^2}_{\text{minimized by}\;
    \bar{x}^* = \bar{b}}
  + \underbrace{\norm{D\tilde\bx}_1 + \frac{1}{2\lambda}
    \norm{\tilde\bx - \tilde{\bb}}_{2}^{2}}_{\text{TV in}\;\tilde\bx}.
\]
Setting $\bz = D\tilde\bx$ and recovering
$\tilde\bx = (D^\top D)^{\dagger}D^\top\bz$ via a Poisson solve
(see Section~\ref{sec:implementation}), the remaining minimization
becomes the LASSO problem
\begin{equation}\label{eq:lasso-final}
  \min_{\bz \in \R^{n-1}} \;
  \norm{\bz}_1 + \frac{1}{2\lambda}\norm{M\bz - \tilde{\bb}}_{2}^{2},
  \qquad M = (D^\top D)^{\dagger}D^\top.
\end{equation}
After solving for $\bz^*$, the full solution is recovered as
$\bx^* = M\bz^* + \bar{b}\bone$.
The difference operator $D \in \R^{(n-1)\times n}$ has rank $n-1$ and is
surjective onto $\R^{n-1}$, so the change of variables $\bz = D\tilde\bx$
is a bijection between $\bone^\perp$ and $\R^{n-1}$ and the LASSO
reduction is exact.

\subsection{Extension to two dimensions}
\label{sec:dim-by-dim}

In 2D, the stacked difference operator
$D = \bigl[\begin{smallmatrix} D_x \\ D_y \end{smallmatrix}\bigr]$
is not surjective, and the LASSO reduction introduces a relaxation
gap.  We instead apply the 1D TV denoiser
dimension-by-dimension: given $U \in \R^{N_x \times N_y}$, first
solve the 1D TV problem for each column, then for each row.
Each sweep preserves the column (resp.\ row) total mass, so the overall
total mass of $U$ is preserved.
The cost is $N_y$ 1D solves of size $N_x$ plus $N_x$ solves
of size $N_y$, each terminating in finitely many breakpoint steps.
A single pass suffices for limiting the 2D solution.

\subsection{Classical iterative solvers}
\label{sec:classical-solvers}

We briefly review two widely used iterative solvers for TV
minimization. Both methods apply to any analysis operator~$D$
(not just the first-difference operator).

\paragraph{Chambolle--Pock accelerated PDHG~\cite{chambolle2011}}
The PDHG algorithm initializes two step sizes $\tau_0 = \sigma_0 = 0.5$ satisfying $\tau_0\sigma_0\norm{D}_{2}^{2} \leq 1$ (since $\norm{D}_2 < 2$) and iterates repeatedly:
\begin{enumerate}[nosep, label=\arabic*.]
  \item \emph{Dual step}: $\by^{k+1}
    = \mathrm{clip}\!\bigl(\by^k + \sigma_k D\hat{\bx}^k,\,-1,\,1\bigr)$
  \item \emph{Primal step}:
    $\bx^{k+1}
    = \dfrac{\bx^k - \tau_k D^\top\by^{k+1} + \tfrac{\tau_k}{\lambda}\bb}
            {1 + \tau_k/\lambda}$
  \item \emph{Acceleration} ($\mu=1/\lambda$):
    $\theta_{k+1} = \tfrac{1}{\sqrt{1+2\tau_k/\lambda}}$,\;
    $\tau_{k+1} = \theta_{k+1}\tau_k$,\;
    $\sigma_{k+1} = \sigma_k/\theta_{k+1}$
  \item \emph{Extrapolation}:
    $\hat{\bx}^{k+1} = \bx^{k+1} + \theta_{k+1}(\bx^{k+1}-\bx^k)$
\end{enumerate}

\paragraph{ADMM (Douglas--Rachford)} The ADMM algorithm splits $\bz = D\bx$ with penalty $\eta = 10/\lambda$ and iterates repeatedly:
\begin{enumerate}[nosep, label=\arabic*.]
  \item \emph{$\bx$-update} (tridiagonal solve):
    $\bigl(\tfrac{1}{\lambda}I + \eta D^\top D\bigr)\bx
    = \tfrac{1}{\lambda}\bb + \eta D^\top(\bz - \mathbf{v})$
  \item \emph{$\bz$-update} (soft-thresholding):
    $\bz = S_{1/\eta}(D\bx + \mathbf{v})$,\;
    where 
    \[
    S_\alpha(w)_j = \sign(w_j)\max(|w_j|-\alpha,\,0)
    \]
  \item \emph{$\mathbf{v}$-update}:
    $\mathbf{v} \leftarrow \mathbf{v} + D\bx - \bz$
\end{enumerate}
Both methods are convergent for the composite minimization of lower
semicontinuous proper convex functions.  Accelerated PDHG attains the
$O(1/k^2)$ rate, and ADMM often converges
faster in practice but at a higher cost per iteration (one tridiagonal
solve instead of two sparse matrix--vector products).  Tuning the ADMM
penalty $\eta$ is crucial, and a poor choice of $\eta$ can slow
convergence by orders of magnitude (see Figure~\ref{fig:tv-spectral} in Section~\ref{sec:tv-spectral}).
ADMM is known as an efficient solver for TV
minimization~\cite{zhang2025}, and a detailed treatment can be found
in~\cite[Section~4.2]{laizhang2026}.

\section{The Langlois--Darbon Differential Inclusion Algorithm}
\label{sec:algorithm}

\subsection{The dual problem and its polyhedral structure}
\label{sec:dual}

The dual problem of~\eqref{eq:lasso-final}, derived by introducing the dual variable
$\bp = (M\bz - \tilde{\bb})/\lambda$ and dualizing the $\ell_1$ term, is
\begin{equation}\label{eq:dual-lasso-tv}
  \min_{\bp \in \R^{n}} \;
  \frac{\lambda}{2}\norm{\bp}_{2}^{2} + \ip{\bp}{\tilde{\bb}}
  \quad \text{subject to} \quad
  \norm{M^\top\bp}_\infty \leq 1.
\end{equation}
The constraint $\norm{M^\top\bp}_\infty \leq 1$ (with $M^\top\bp \in \R^{n-1}$) defines a
closed convex polyhedron, the same polyhedral structure as in the standard LASSO.

\subsection{Implementation}
\label{sec:implementation}
The DI algorithm solves the dual problem~\eqref{eq:dual-lasso-tv} by computing the \emph{breakpoints} of the piecewise continuous trajectory of its associated differential inclusions. At a feasible point $\bp$, the \emph{equicorrelation set} $\mathcal{E}(\bp) = \{j : \abs{(M^\top\bp)_j} = 1\}$ identifies the active constraints along which the trajectory evolves. The algorithm proceeds as follows:
\begin{enumerate}[label=\textbf{Step~\arabic*.}, leftmargin=3.5em]
  \item \textbf{Offline precomputation.}
    Form the dense matrix $M = D^\dagger = (D^\top D)^\dagger D^\top
    \in \R^{n \times (n-1)}$, which is independent of the data $\bb$.

  \item \textbf{Initialization.}
    Set $\tilde{\bb} = \bb - \operatorname{mean}(\bb)\cdot\bone$.
    If $\tilde{\bb} = \bzero$ (constant input), the solution is
    $\bx^* = \bar{b}\,\bone$ and the algorithm terminates.
    Otherwise, set $\bp_0 = -\tilde{\bb}/\norm{M^\top\tilde{\bb}}_\infty$
    (the unique scaling that places $\bp_0$ on the boundary of the polyhedron)
    and compute $\mathcal{E}_0$.

  \item \textbf{NNLS subproblem.}
    At breakpoint $\bp_k$ with active set $\mathcal{E}_k$
    and signs $s_j = \sign(-(M^\top\bp_k)_j)$,
    form the columns $\tilde{M} = [s_j\, M(:,j)]_{j \in \mathcal{E}_k}$
    and solve the NNLS problem
    $\min_{\bu \geq 0}\,\norm{\tilde{M}\bu - (\lambda\bp_k + \tilde{\bb})}_{2}^{2}$.
    The descent direction is
    $\bd_k = \tilde{M}\bu - (\lambda\bp_k + \tilde{\bb})$.

  \item \textbf{Ratio test.}
    Compute the step $\Delta_*$ to the next breakpoint:
    the smallest $\Delta > 0$ at which some inactive constraint
    $|(M^\top(\bp_k + \Delta\bd_k))_j| = 1$ becomes active.
    If $\lambda\Delta_* \geq 1$, the algorithm has converged.
    Otherwise, update $\bp_{k+1} = \bp_k + \Delta_*\,\bd_k$.

  \item \textbf{Recovery.}
    Set $z^*_j = s_j u_j$ for $j \in \mathcal{E}$
    (zero elsewhere) and recover
    $\bx^* = M\bz^* + \operatorname{mean}(\bb)\cdot\bone$.
\end{enumerate}

\noindent
Finite termination is guaranteed: $\norm{\bd_k}_{2}$ strictly decreases along each face of the polyhedron defined by $\norm{M^{\top}\bp}_{\infty} \leq 1$, so the algorithm
terminates~\cite{langlois2025}.

\paragraph{Forming $M$ via the eigendecomposition of $D^\top D$}
The matrix $D^\top D$ is the 1D Neumann Laplacian, a tridiagonal
$n \times n$ matrix whose eigenvalues and orthonormal eigenvectors are known in
closed form \cite{laizhang2026}:
\[
  \mu_k = 2 - 2\cos\!\bigl(\pi k/n\bigr), \qquad
  v_k(j) = c_k \cos\!\Bigl(\frac{k\pi(j - \tfrac{1}{2})}{n}\Bigr),
\]
for $j = 1, \ldots, n$ and $k = 0, \ldots, n-1$,
where $c_0 = 1/\sqrt{n}$ and $c_k = \sqrt{2/n}$ for $k \geq 1$.
The eigenvectors correspond to the DCT-II basis on the one-half grid
$x_j = (j - \tfrac{1}{2})/n$.
Setting $1/\mu_0 = 0$ (pseudoinverse convention for the zero
eigenvalue), we compute
\[
  M = (D^\top D)^\dagger D^\top
  = C^\top \diag(\mu_0^{-1}, \ldots, \mu_{n-1}^{-1})\, C\, D^\top,
\]
where $C$ is the orthogonal matrix with entries $C_{kj} = c_k \cos\bigl(k\pi(j - \tfrac{1}{2})/n\bigr)$.
This offline computation costs $O(n^2\log n)$ via the
DCT (or $O(n^3)$ with dense matrix multiplications) and produces a dense
$n \times (n-1)$ matrix (less than 8\,MB for $n = 1000$).
Once $M$ is stored, column extraction
$M(:,\mathcal{E})$ is a trivial indexing operation, and the
matrix--vector products $M\bz$ and $M^\top\bp$ needed at each
breakpoint are standard dense matvecs.

\paragraph{Warm-started NNLS via Meyer's algorithm}
The NNLS subproblem computed at each breakpoint is the computational
bottleneck.  We solve it using Meyer's
algorithm~\cite{meyer2013simple}, a generalization of the classical
Lawson--Hanson active-set method that accepts an arbitrary starting
active set.  Because the equicorrelation set $\mathcal{E}$ typically
grows by only one index per breakpoint, the passive set from the
previous solve provides a near-optimal warm start: the NNLS needs
only $O(1)$ active-set pivots instead of $O(|\mathcal{E}|)$ from a
cold start.  In addition, columns of $M$ that have been extracted in
previous breakpoints are cached so that only newly activated indices
require a column lookup.  The combination yields a roughly $5\times$
speedup at $n = 1000$.

\paragraph{Complexity}
Each breakpoint step involves one NNLS solve and one ratio test ($O(n^2)$ for computing $M^\top\bd$ and checking constraints). A cold-start NNLS solve costs on the order of $|\mathcal{E}|^2 n$ flops with active-set pivots, but with the warm start described above the NNLS typically needs only $O(1)$ pivots of $O(|\mathcal{E}|\,n)$ flops each, so the observed cost per breakpoint step is $O(n^2)$. Over $K$ breakpoint steps the total empirical cost is therefore $O(K n^2)$, and since $K = O(n)$ in all our tests ($K \le 1.5n$; see Tables~\ref{tab:di-summary} and~\ref{tab:fd4-summary}), the empirical worst-case complexity is $O(n^3)$. Although the worst-case complexity of the NNLS problem with active-set methods remains open, in practice, for clean piecewise-constant data such as the truncated Fourier sums of step functions in Section~\ref{sec:tv-spectral}, the number of breakpoint steps is much smaller than $n$ and depends on signal complexity rather than grid size (see Table~\ref{tab:di-summary}).

\subsection{Numerical comparison}
\label{sec:denoising}

We compare the DI algorithm with the iterative methods of
Section~\ref{sec:classical-solvers} on TV denoising problems.
All three methods minimize the same 1D TV objective
$F(\bx) = \norm{D\bx}_1 + \frac{1}{2\lambda}\norm{\bx-\bb}_{2}^{2}$.
The DI algorithm is implemented through the equivalent LASSO problem $\min_{\bz \in \R^{n-1}} \norm{\bz}_1 + \frac{1}{2\lambda}\norm{M\bz - \tilde{\bb}}_{2}^{2}$
with $M = D^\dagger$ and $\tilde{\bb} = \bb - \bar{b}\,\bone$.
Since $D$ is surjective, this reformulation is exact.

\subsubsection{Gaussian noise denoising}
\label{sec:tv-gaussian}

We first consider a piecewise constant signal corrupted by additive
Gaussian noise. The true signal on $n = 200$ grid points has five
constant segments with values $\{1, -0.5, 2, 0, -1\}$ and four jumps.
The observation is $\bb = \bx_{\mathrm{true}} + \boldsymbol{\varepsilon}$,
where $\boldsymbol{\varepsilon} \sim \mathcal{N}(0, 0.3^2 I)$, and we
use $\lambda = 1$.

\begin{figure}[!htbp]
  \centering
  \begin{subfigure}[t]{0.48\textwidth}
    \centering
    \includegraphics[width=\textwidth]{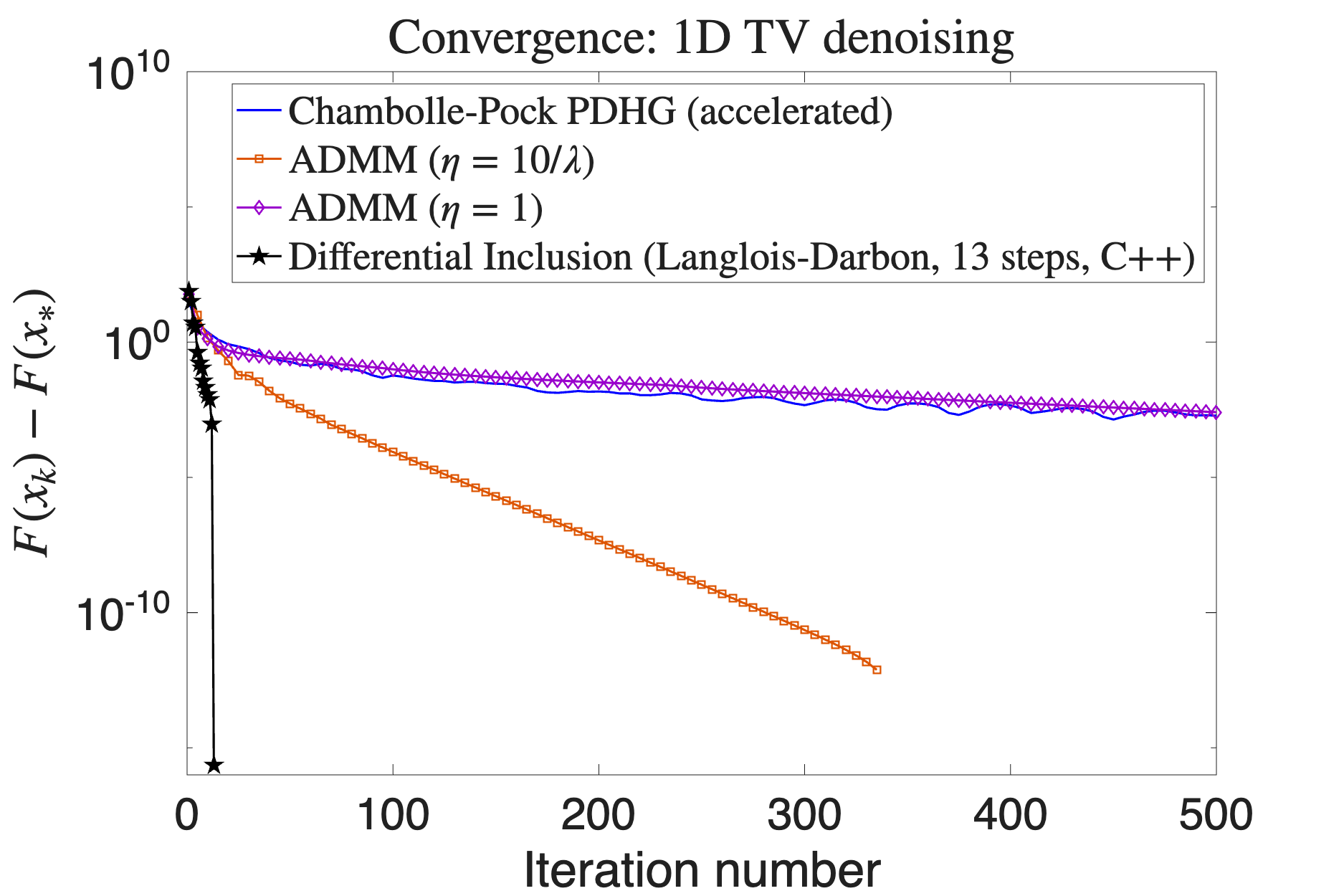}
    \caption{Objective error vs.\ iteration number.}
    \label{fig:tv-gaussian-iter}
  \end{subfigure}
  \hfill
  \begin{subfigure}[t]{0.48\textwidth}
    \centering
    \includegraphics[width=\textwidth]{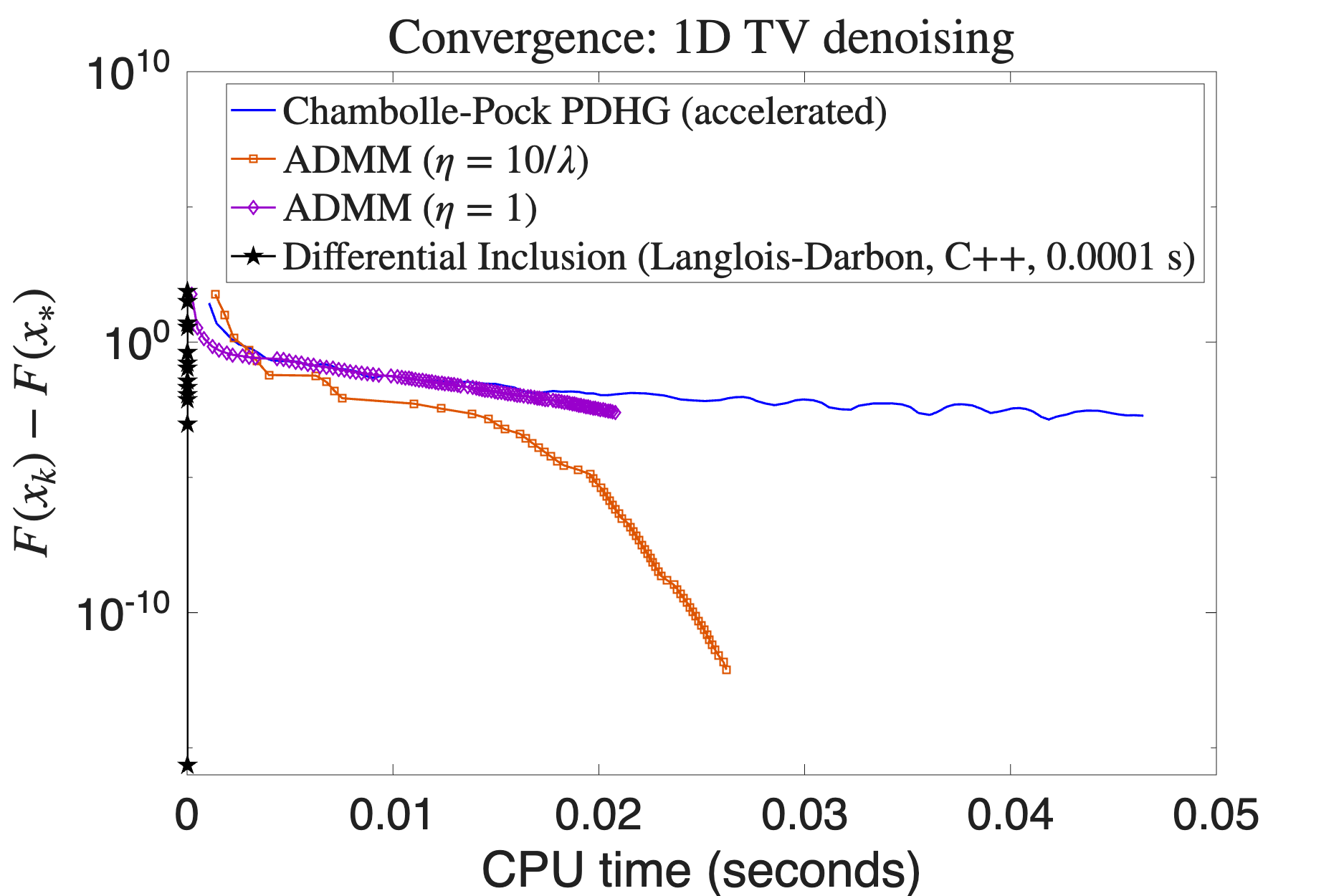}
    \caption{Objective error vs.\ CPU time.}
    \label{fig:tv-gaussian-time}
  \end{subfigure}
  \caption{Convergence comparison for TV denoising of a piecewise
  constant signal with Gaussian noise ($n=200$, $\sigma=0.3$, $\lambda=1$).
  The DI algorithm converges in 13~steps.}
  \label{fig:tv-gaussian}
\end{figure}

The DI algorithm converges in
13~breakpoint steps.
As shown in Figure~\ref{fig:tv-gaussian}, ADMM with tuned
$\eta = 10/\lambda$ reaches ${\sim}\,10^{-13}$ within 500 iterations,
while PDHG and untuned ADMM stall near $10^{-3}$. At $n = 1000$ with
the same noise level and proportionally scaled jump locations, the
solver converges in 59~steps.

\subsubsection{Spectral Gibbs denoising: piecewise constant signal}
\label{sec:tv-spectral}

Consider removing Gibbs oscillations from a truncated
Fourier series.
The true signal is a piecewise constant function on $[0,1)$
with two jumps:
\[
  f(x) = \begin{cases}
    +1 & x \in [0, 0.25) \cup [0.75, 1),\\
    -1 & x \in [0.25, 0.75).
  \end{cases}
\]
The ``noisy'' signal $\bb$ is the $N$-term partial Fourier sum of $f$,
evaluated at $n = N$ uniform grid points $x_j = j/N$:
\[
  b_j = \sum_{m=-N/2+1}^{N/2} \hat{f}_m \, e^{2\pi i m x_j},
  \qquad \hat{f}_m = \int_0^1 f(x)\, e^{-2\pi i m x}\, dx.
\]
(The real part of the partial sum is taken, so that the unpaired
mode $m = N/2$ does not introduce a spurious imaginary component.)
This is the output of a Fourier spectral method on an
$N$-point grid, with Gibbs oscillations near each discontinuity
($\norm{\bb - \bx_{\mathrm{true}}}_2 \approx 1.49$).

We apply the TV denoising model
$F(\bx) = \norm{D\bx}_1 + \frac{1}{2\lambda}\norm{\bx-\bb}_{2}^{2}$
with $\lambda = 0.3$. This is smaller than the choice $\lambda = 1$
used in Section~\ref{sec:tv-gaussian}, reflecting the lower effective
noise level. Figure~\ref{fig:tv-spectral-signal}
(Section~\ref{sec:formulation}) shows that TV denoising removes the
Gibbs oscillations and recovers the piecewise constant structure.

\begin{figure}[!htbp]
  \centering
  \begin{subfigure}[t]{0.48\textwidth}
    \centering
    \includegraphics[width=\textwidth]{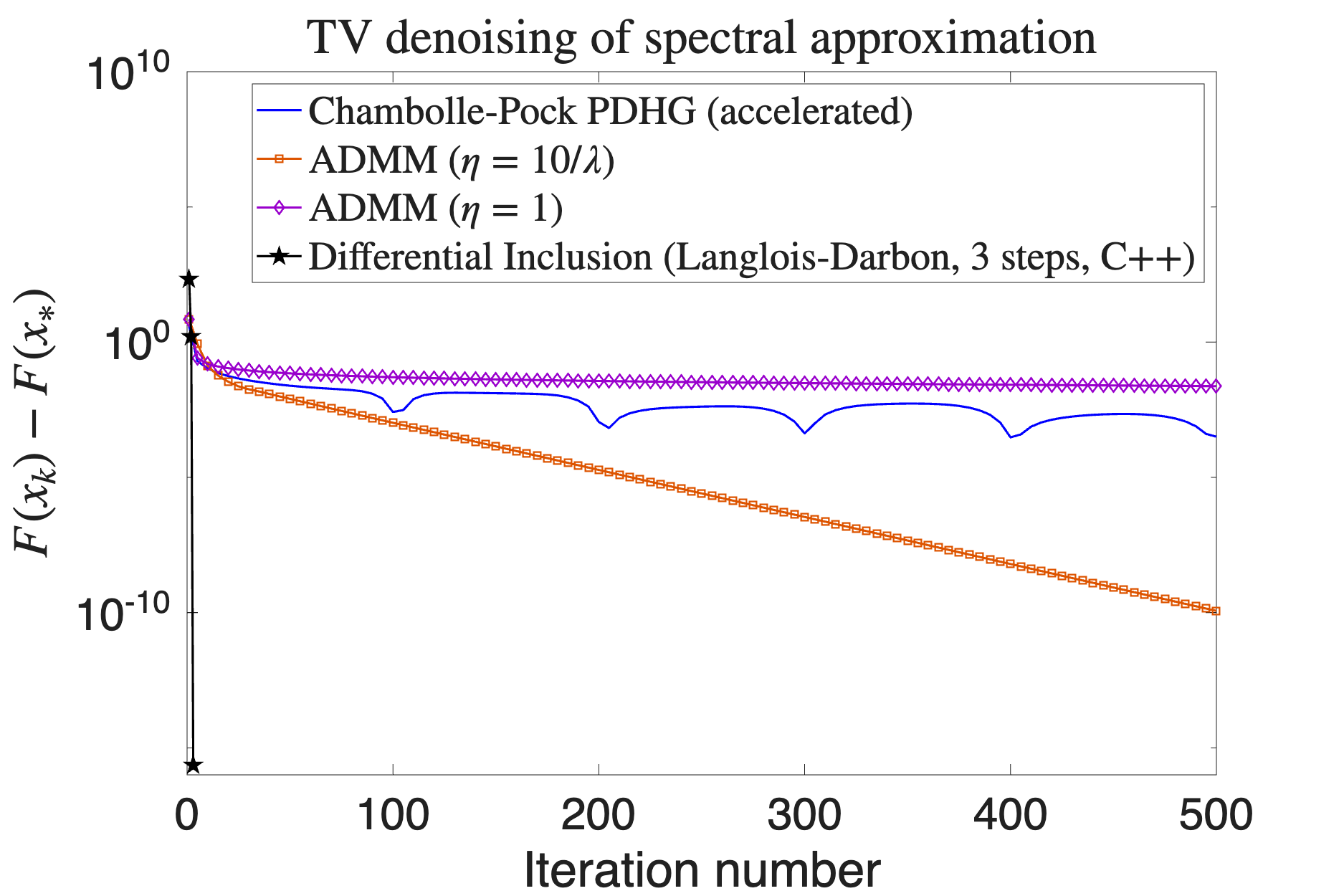}
    \caption{$N=200$: objective error vs.\ iteration.}
    \label{fig:tv-spectral-iter}
  \end{subfigure}
  \hfill
  \begin{subfigure}[t]{0.48\textwidth}
    \centering
    \includegraphics[width=\textwidth]{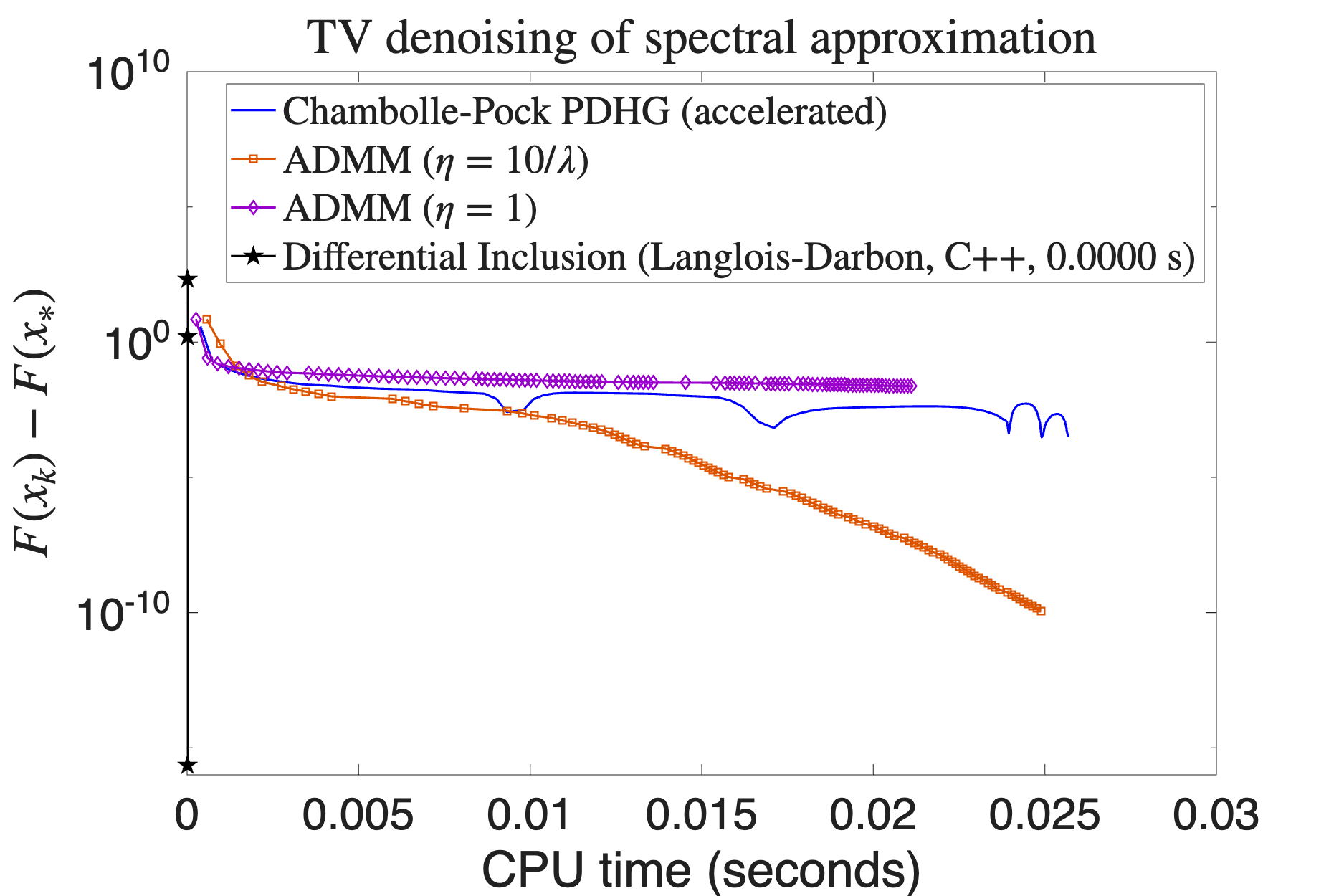}
    \caption{$N=200$: objective error vs.\ CPU time.}
    \label{fig:tv-spectral-time}
  \end{subfigure}
  \\[6pt]
  \begin{subfigure}[t]{0.48\textwidth}
    \centering
    \includegraphics[width=\textwidth]{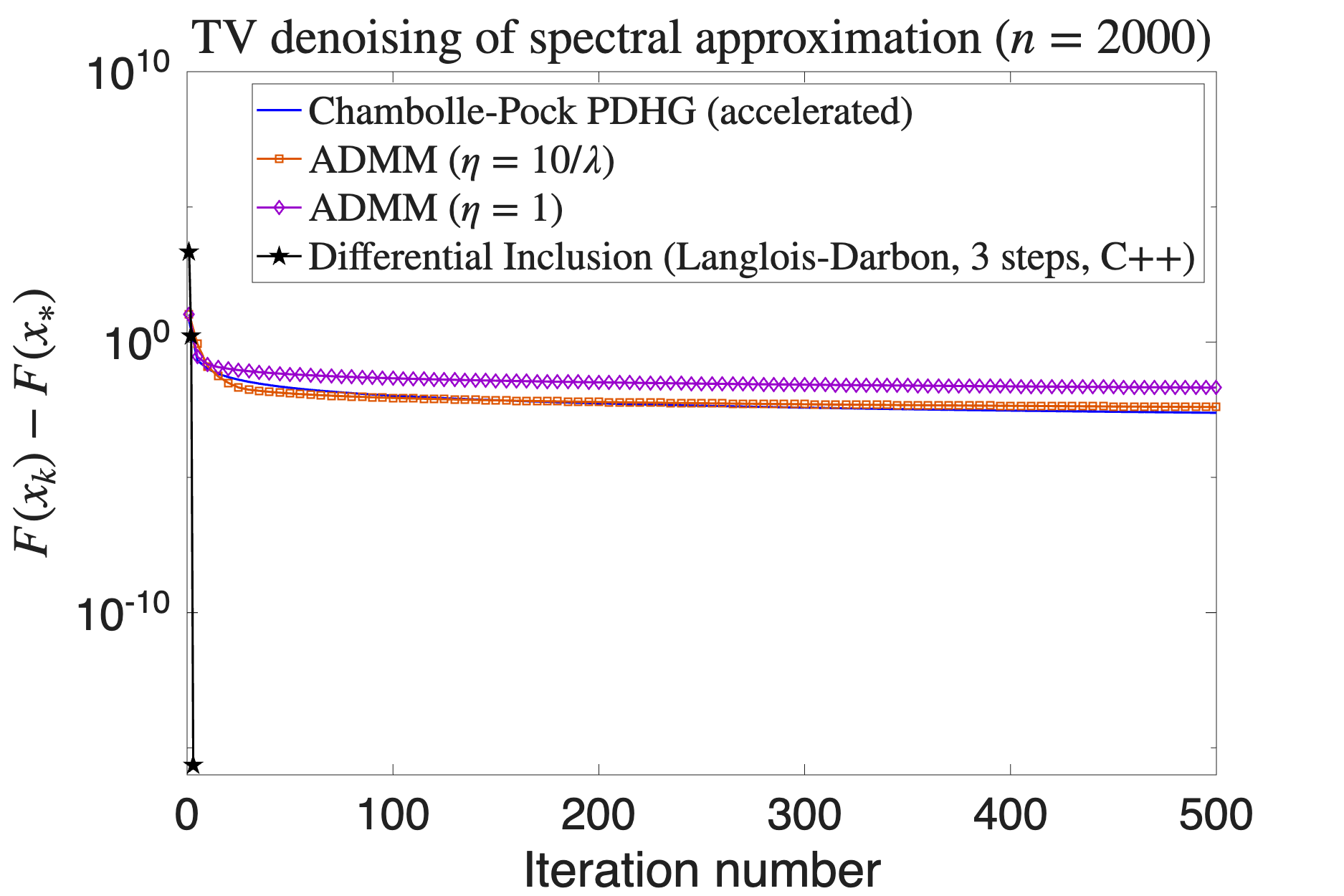}
    \caption{$N=2000$: objective error vs.\ iteration.}
    \label{fig:tv-spectral-large-iter}
  \end{subfigure}
  \hfill
  \begin{subfigure}[t]{0.48\textwidth}
    \centering
    \includegraphics[width=\textwidth]{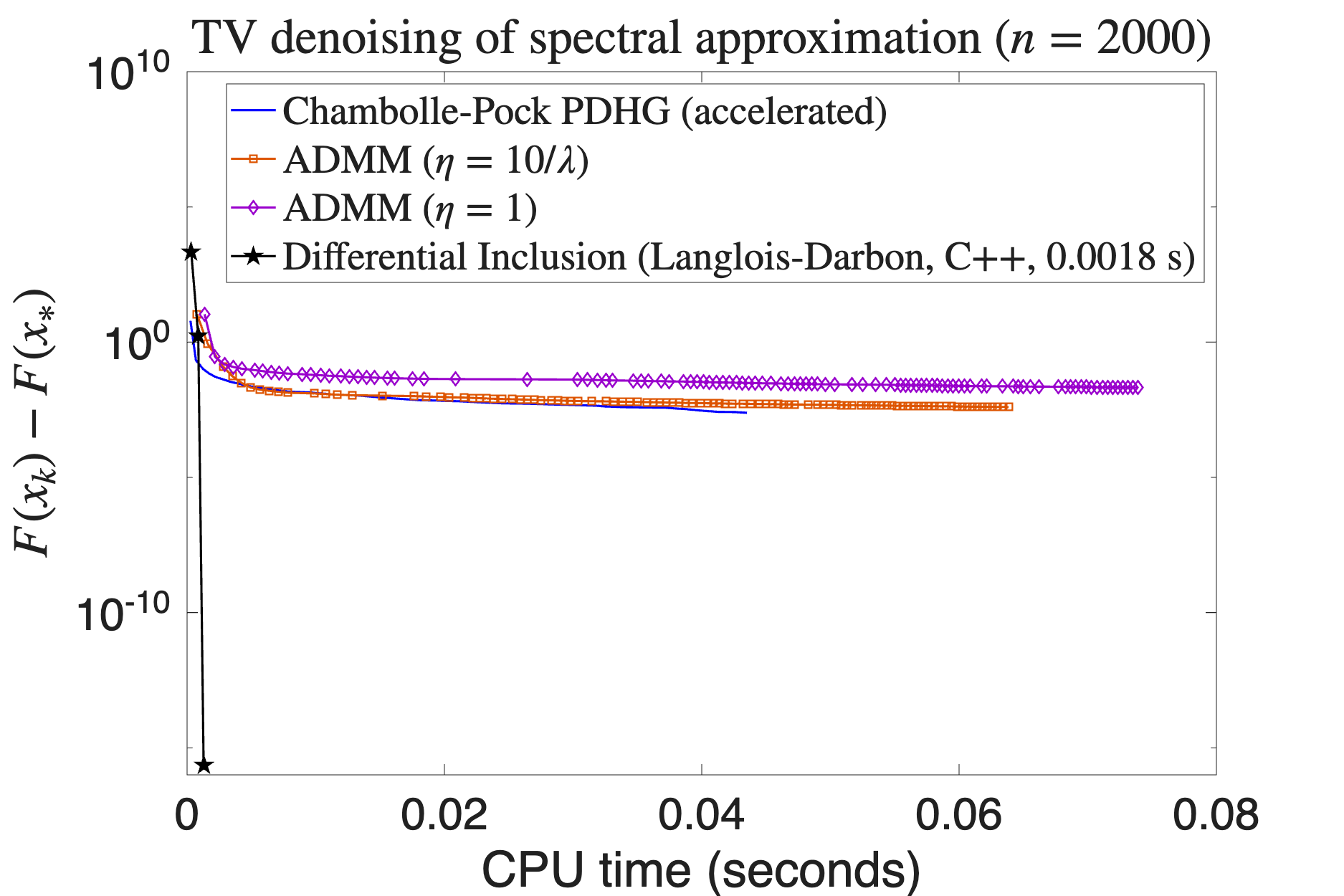}
    \caption{$N=2000$: objective error vs.\ CPU time.}
    \label{fig:tv-spectral-large-time}
  \end{subfigure}
  \caption{Convergence comparison for TV denoising of the
  partial Fourier sum of a piecewise constant signal ($\lambda=0.3$).
  The DI algorithm converges in 3 steps at both
  $N=200$ and $N=2000$.}
  \label{fig:tv-spectral}
\end{figure}

At $N = 200$, the DI algorithm converges in
3 breakpoint steps. At $N = 2000$, it again converges in
3 breakpoint steps. The step count does not increase with
$N$ because the underlying piecewise constant structure is unchanged.
The iterative methods stall near $10^{-2}$--$10^{-3}$ within
500 iterations (see Figure~\ref{fig:tv-spectral}).

\FloatBarrier
\subsubsection{Spectral Gibbs denoising: piecewise smooth signal}
\label{sec:tv-spectral2}

The previous test used a piecewise \emph{constant} signal, for which
the DI algorithm required only a few breakpoint steps.
Consider a piecewise \emph{smooth} signal:
\[
  f(x) = \begin{cases}
    \phantom{-}\sin(\pi x) & x \in [0, 0.5),\\
    -\sin(\pi x) & x \in [0.5, 1),
  \end{cases}
\]
which has a jump discontinuity of size~2 at $x = 0.5$
(from $\sin(\pi/2) = 1$ to $-\sin(\pi/2) = -1$),
while the smooth portions are sinusoidal arcs.
As before, the ``noisy'' signal $\bb$ is the $N$-term partial
Fourier sum of $f$ evaluated at $n = N$ uniform grid points,
producing Gibbs oscillations near $x=0.5$.
We apply TV denoising with $\lambda = 1$.

\begin{figure}[!htbp]
  \centering
  \begin{subfigure}[t]{0.48\textwidth}
    \centering
    \includegraphics[width=\textwidth]{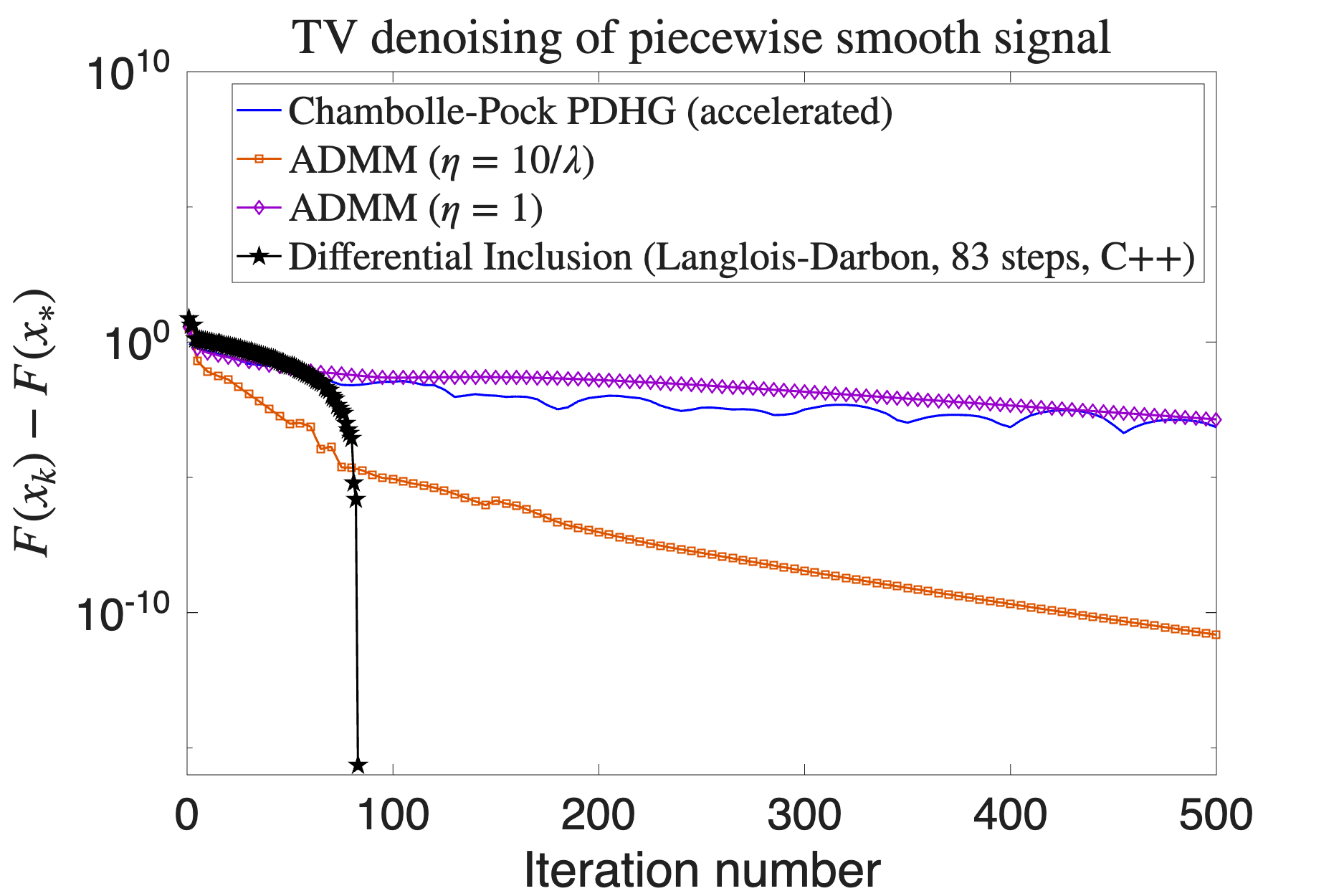}
    \caption{$N=200$: objective error vs.\ iteration.}
    \label{fig:tv-spectral2-iter}
  \end{subfigure}
  \hfill
  \begin{subfigure}[t]{0.48\textwidth}
    \centering
    \includegraphics[width=\textwidth]{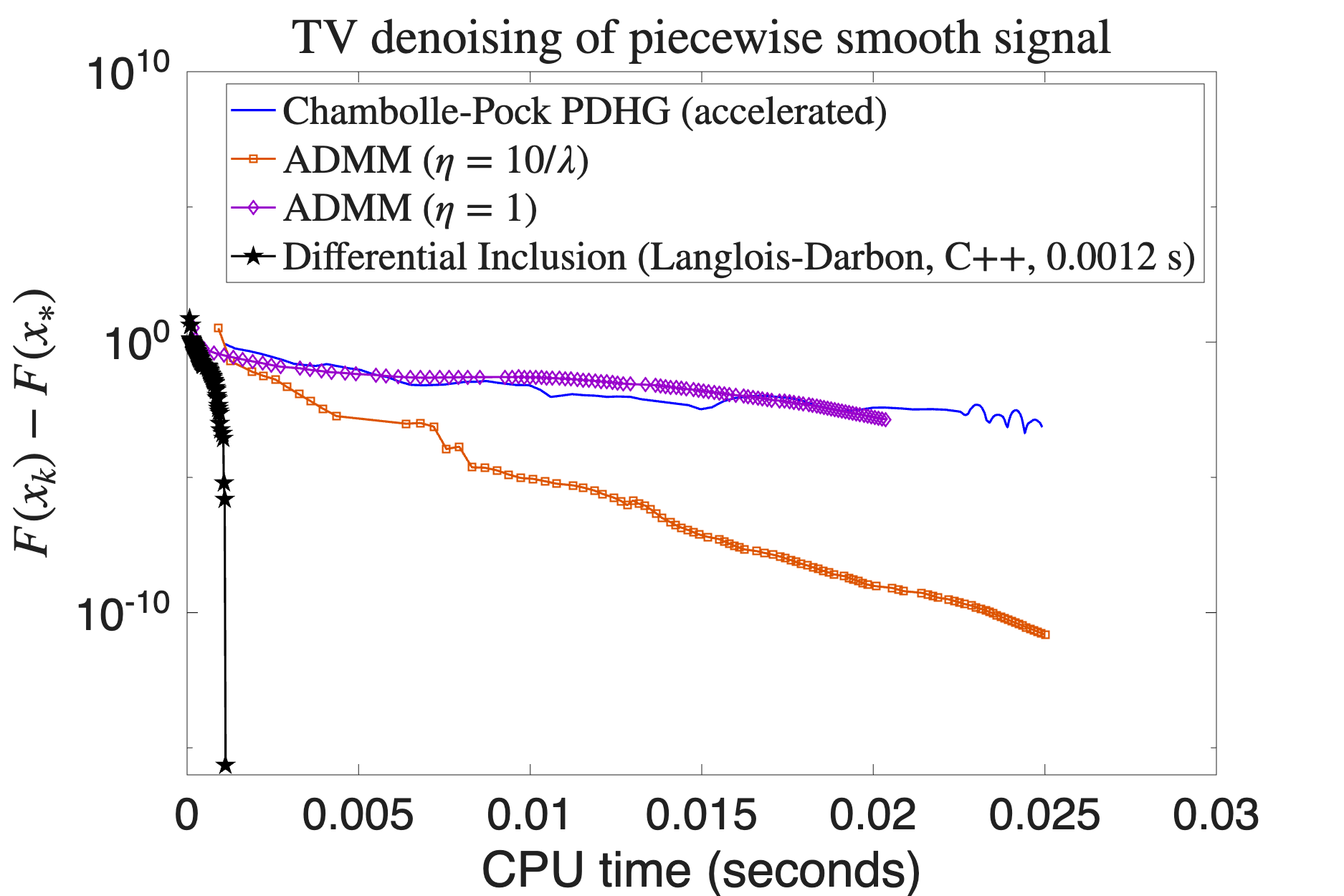}
    \caption{$N=200$: objective error vs.\ CPU time.}
    \label{fig:tv-spectral2-time}
  \end{subfigure}
  \\[6pt]
  \begin{subfigure}[t]{0.48\textwidth}
    \centering
    \includegraphics[width=\textwidth]{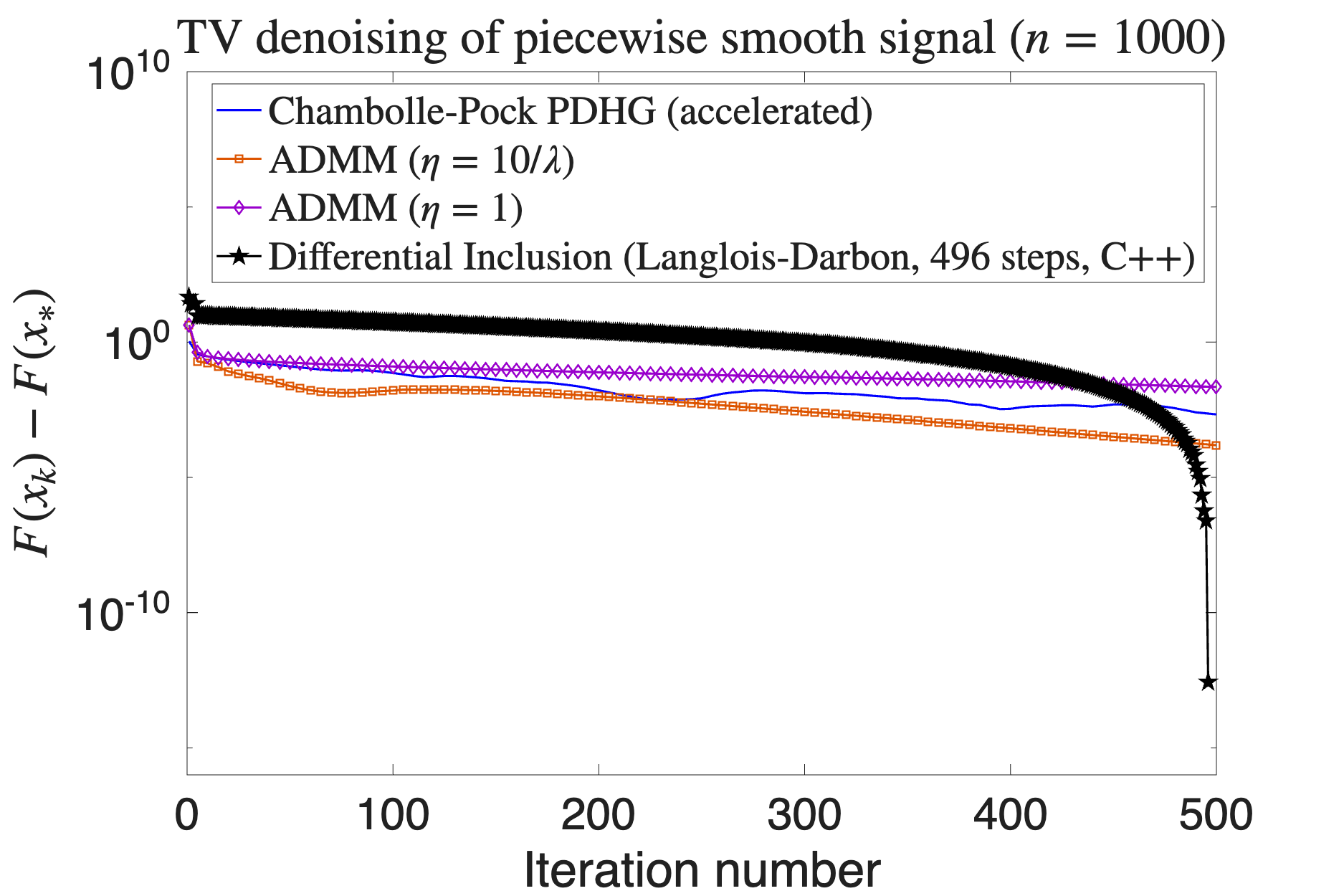}
    \caption{$N=1000$: objective error vs.\ iteration.}
    \label{fig:tv-spectral2-large-iter}
  \end{subfigure}
  \hfill
  \begin{subfigure}[t]{0.48\textwidth}
    \centering
    \includegraphics[width=\textwidth]{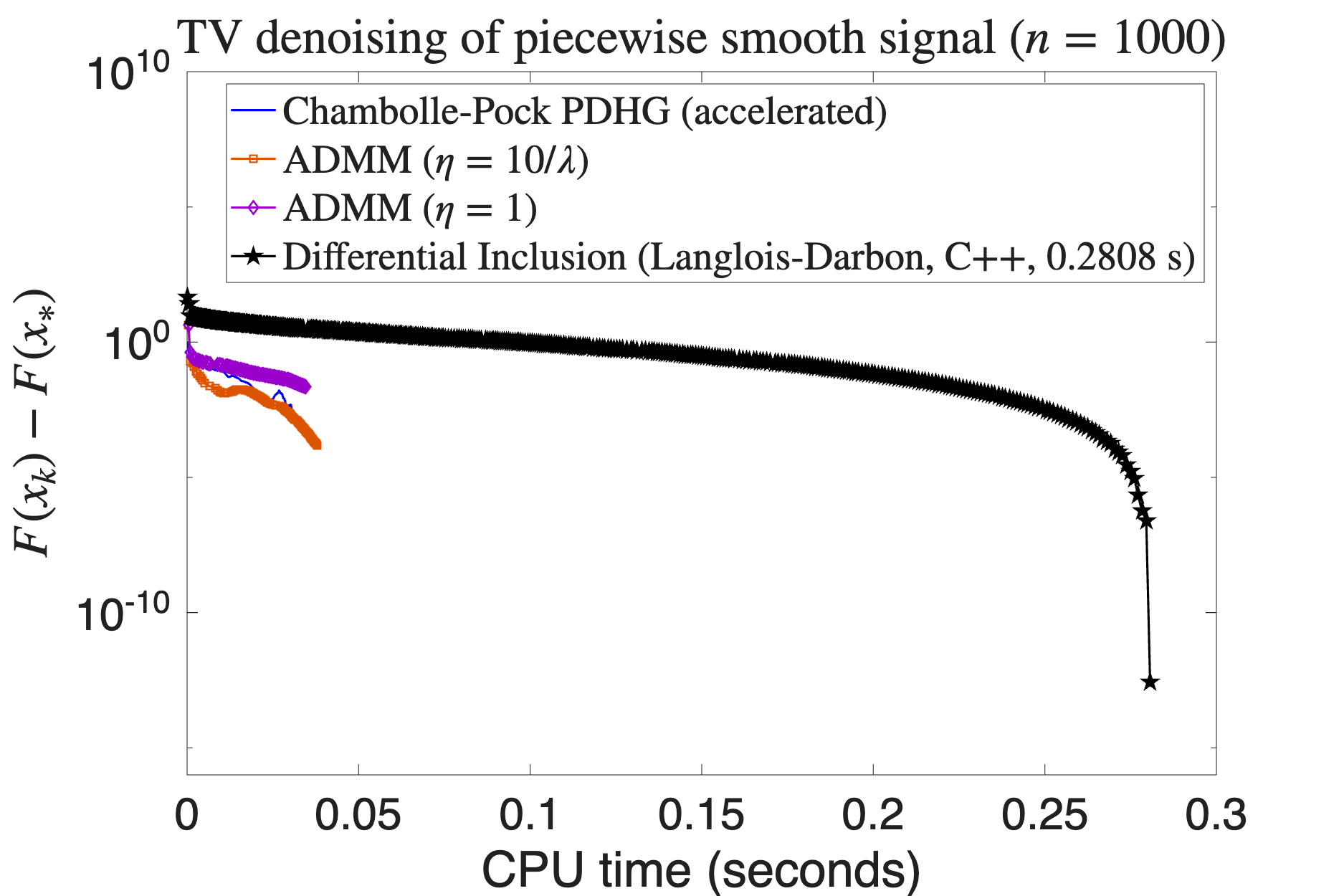}
    \caption{$N=1000$: objective error vs.\ CPU time.}
    \label{fig:tv-spectral2-large-time}
  \end{subfigure}
  \caption{Convergence comparison for TV denoising of the
  partial Fourier sum of a piecewise smooth signal ($\lambda=1$).
  }
  \label{fig:tv-spectral2}
\end{figure}

As shown in Figure~\ref{fig:tv-spectral2-signal}
(Section~\ref{sec:formulation}), the Gibbs oscillations are eliminated
while the sinusoidal profile is retained in the smooth
regions.  That figure uses a milder regularization
($\lambda = 0.08$) chosen for visual clarity; the $\lambda = 0.08$ setting is the one timed in Table~\ref{tab:di-summary}, while the convergence study
here uses $\lambda = 1$ (Figure~\ref{fig:tv-spectral2}).
At $N = 200$, the DI algorithm requires
83~breakpoint steps, compared with 3 steps for the
piecewise constant signal. The TV-optimal reconstruction of a smooth
signal has many nonzero differences $D\bx^*$, so the equicorrelation
set grows to a large fraction of the $n - 1 = 199$ possible indices.
For a piecewise constant signal with $k$ jumps, the equicorrelation set
typically has size~${\sim}\,k$. PDHG and untuned ADMM stall near
$10^{-3}$ after 500 iterations.

At $N = 1000$ with $\lambda = 1$, the DI algorithm
requires 496~breakpoint steps. The increase reflects
the finer grid, which produces more nonzero differences in the
TV-optimal solution and enlarges the equicorrelation set. ADMM with
tuned $\eta = 10/\lambda$ reaches ${\sim}\,10^{-6}$, while PDHG and
untuned ADMM stall near $10^{-2}$. Figure~\ref{fig:tv-spectral2}
shows the convergence for both $N = 200$ and $N = 1000$.

At first glance, Figures~\ref{fig:tv-spectral2-time}
and~\ref{fig:tv-spectral2-large-time} may suggest that ADMM and PDHG
are competitive in CPU time. However,
Figure~\ref{fig:tv-spectral2-large-iter} shows that at $N = 1000$ both
iterative methods remain far from convergence after 500 iterations.
Each step of the DI algorithm is more expensive
than one ADMM or PDHG iteration, since it involves an NNLS solve.
Even so, ADMM and PDHG need substantially more iterations to
reach the same accuracy.

\section{TV Limiter for Conservation Laws}
\label{sec:pde}

We use TV denoising as a \emph{limiter} in time-dependent PDE
computations. We consider two classes of base schemes: a fifth-order
upwind finite difference scheme (FD5), used for scalar Burgers'
equation (Section~\ref{sec:pde-burgers}), the 1D compressible Euler
equations (Section~\ref{sec:fd-euler}), and a 2D Riemann problem for
compressible Euler (Section~\ref{sec:pde-riemann2d}); and Fourier
pseudospectral methods for two-dimensional scalar convection and
incompressible Euler (Sections~\ref{sec:pde-rotation}--\ref{sec:pde-euler}).
When applied as a per-step limiter, the TV denoising
step is invoked every $K$ time steps, but only when the current
total variation $\mathrm{TV}(u^n)$ exceeds a threshold
$\mathrm{TV}_{\mathrm{thresh}} = 1.2\,\mathrm{TV}(u_0)$,
where $\mathrm{TV}(u_0)$ is the total variation of the initial condition.

\FloatBarrier
\subsection{FD5 for Burgers' equation: shock capturing}
\label{sec:pde-burgers}

Consider Burgers' equation with periodic boundary conditions:
\[
  u_t + \Bigl(\frac{u^2}{2}\Bigr)_x = 0, \qquad x \in [0,2\pi),
  \quad u(x,0) = 1 + \sin(x).
\]
The shock forms at
$t_s = -1/\min_x u_0'(x) = 1$ and propagates at speed~1.
We discretize with $N = 256$ uniform grid points
$x_i = i\,\Delta x$, $\Delta x = 2\pi/N$, using the fifth-order
linear upwind conservative finite difference scheme (FD5)
\begin{equation}\label{eq:fd5}
  \frac{d}{dt} u_i
  = -\frac{1}{\Delta x}
    \bigl(\hat{f}_{i+1/2} - \hat{f}_{i-1/2}\bigr),
\end{equation}
with Lax--Friedrichs flux splitting~\cite{Shu1998,JiangShu1996}.
The flux is split as $f^{\pm}(u) = \tfrac{1}{2}(f(u) \pm \alpha u)$
with $\alpha = \max|f'(u)|$, and the numerical flux
$\hat{f}_{i+1/2} = \hat{f}^{+}_{i+1/2} + \hat{f}^{-}_{i+1/2}$
is reconstructed as
\[
  \hat{f}^{+}_{i+1/2}
  = \frac{1}{60}\bigl(2f^{+}_{i-2} - 13f^{+}_{i-1}
    + 47f^{+}_{i} + 27f^{+}_{i+1} - 3f^{+}_{i+2}\bigr),
\]
\[
  \hat{f}^{-}_{i+1/2}
  = \frac{1}{60}\bigl(-3f^{-}_{i-1} + 27f^{-}_{i}
    + 47f^{-}_{i+1} - 13f^{-}_{i+2} + 2f^{-}_{i+3}\bigr).
\]
This is the linear scheme to
which WENO5~\cite{JiangShu1996} reduces when the optimal
linear weights are used in place of the nonlinear WENO weights.
We integrate~\eqref{eq:fd5} in time with
RK4  (time step
$\Delta t = 0.5\,\Delta x$, so the CFL number $\Delta t\,\max|u|/\Delta x
\approx 1$ since $\max|u| \approx 2$; 123~steps to $T = 1.5$).

For this test case, the FD5 scheme is stable without any limiter:
the no-limiter solution has $L^2$ error $5.75 \times 10^{-2}$ and
$\mathrm{TV} = 5.83$ (versus exact $\mathrm{TV} = 3.99$), with
oscillations confined near the shock.  TV post-processing with
different regularization strengths reduces the error:
\begin{center}
\begin{tabular}{lcc}
  \hline
  Method & $L^2$ error & TV \\
  \hline
  FD5, no limiter         & $5.75 \times 10^{-2}$ & 5.83 \\
  TV post-proc $\lambda{=}0.02$ & $5.53 \times 10^{-2}$ & 5.23 \\
  TV post-proc $\lambda{=}0.05$ & $5.31 \times 10^{-2}$ & 4.76 \\
  TV post-proc $\lambda{=}0.10$ & $5.11 \times 10^{-2}$ & 4.29 \\
  \hline
\end{tabular}
\end{center}
At $\lambda = 0.10$, TV post-processing reduces the $L^2$ error by 11\%
and brings the total variation closer to the exact value.
The total mass is preserved to machine precision
($|\bar{u}_h - \bar{u}_0| < 10^{-15}$).

\begin{figure}[!htbp]
  \centering
  \begin{subfigure}[t]{0.48\textwidth}
    \centering
    \includegraphics[width=\textwidth]{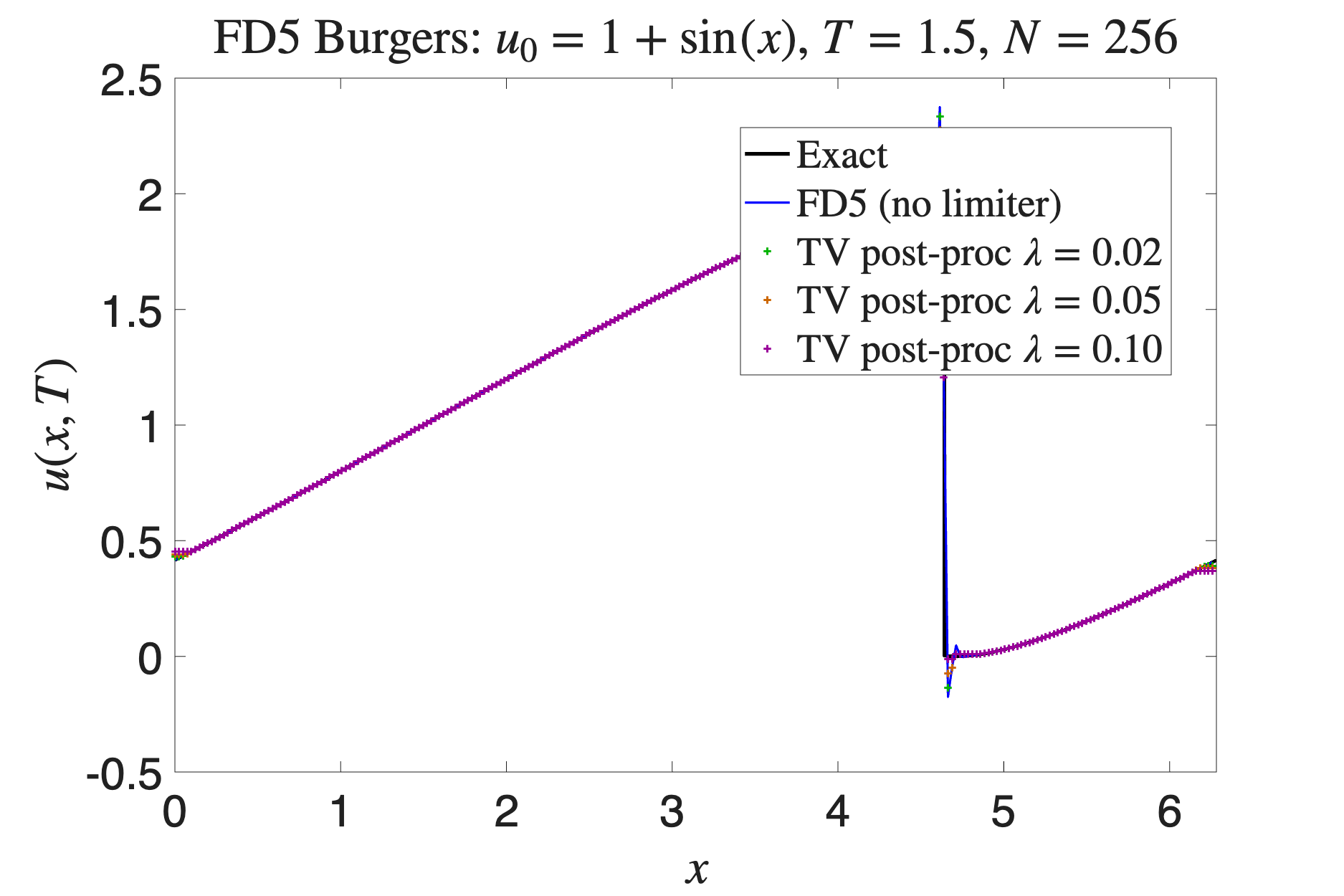}
    \caption{Solution at $T = 1.5$.}
    \label{fig:burgers-shock-sol}
  \end{subfigure}
  \hfill
  \begin{subfigure}[t]{0.48\textwidth}
    \centering
    \includegraphics[width=\textwidth]{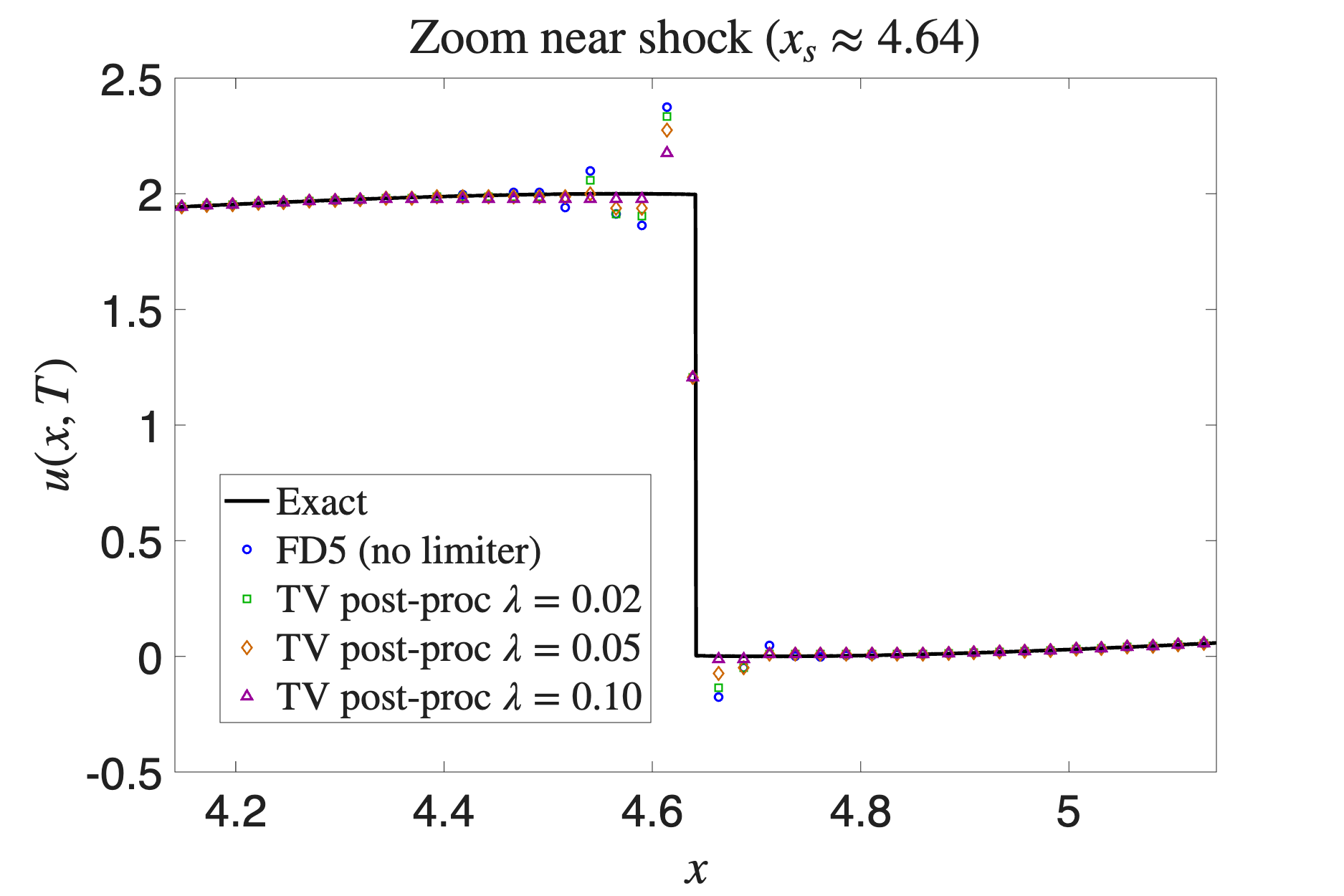}
    \caption{Zoom near the shock.}
    \label{fig:burgers-shock-zoom}
  \end{subfigure}
  \caption{FD5 Burgers' equation with $u_0 = 1 + \sin(x)$
  at $T = 1.5$ ($N = 256$).
  Blue: FD5 without post-processing.
  Colored markers: TV post-processing with $\lambda = 0.02, 0.05, 0.10$.
  The FD5 scheme is stable for this test case. TV post-processing removes
  the remaining oscillations near the shock.}
  \label{fig:burgers-shock}
\end{figure}

\begin{figure}[!htbp]
  \centering
  \begin{subfigure}[t]{0.48\textwidth}
    \centering
    \includegraphics[width=\textwidth]{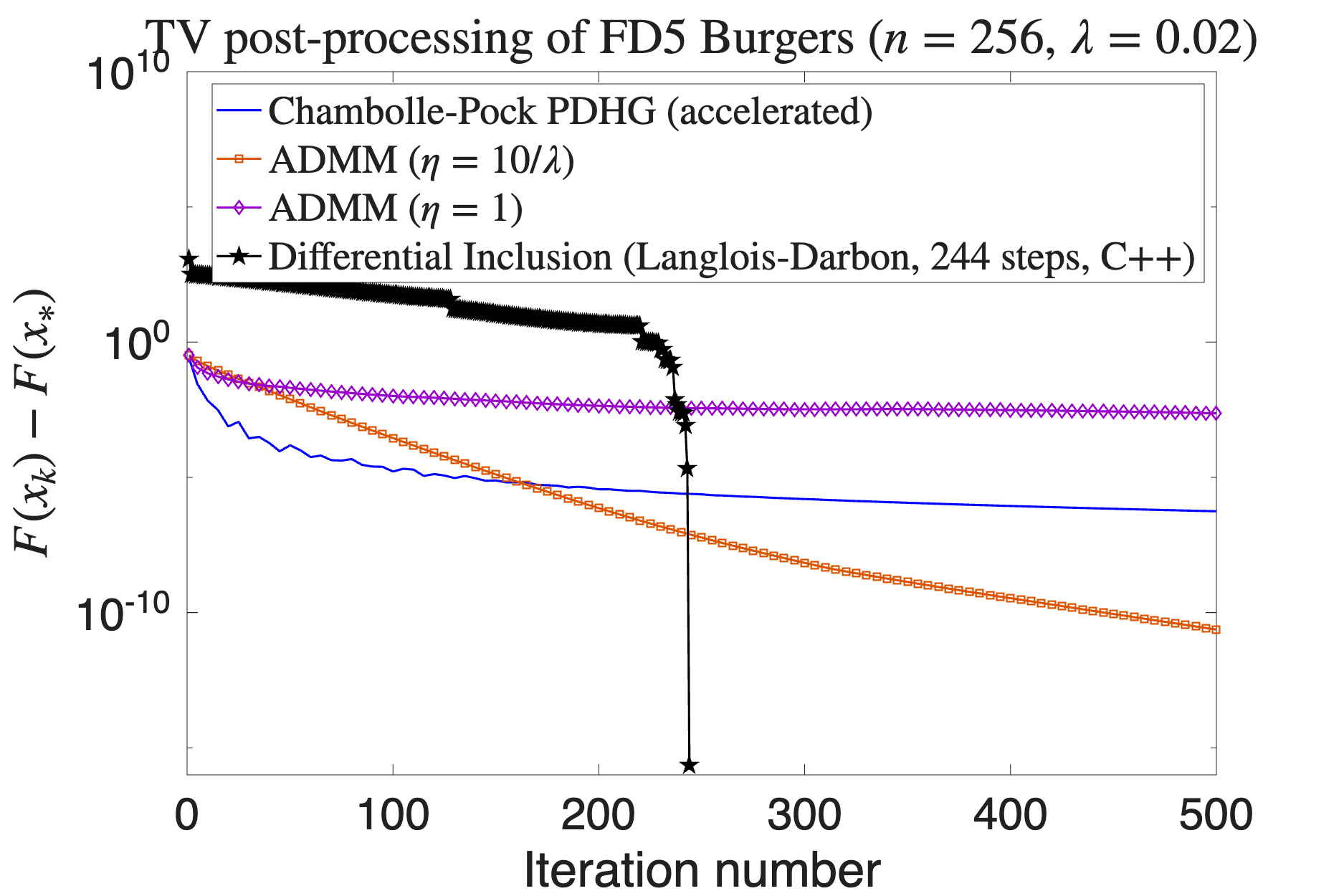}
    \caption{Objective error vs.\ iteration number.}
    \label{fig:burgers-shock-conv-iter}
  \end{subfigure}
  \hfill
  \begin{subfigure}[t]{0.48\textwidth}
    \centering
    \includegraphics[width=\textwidth]{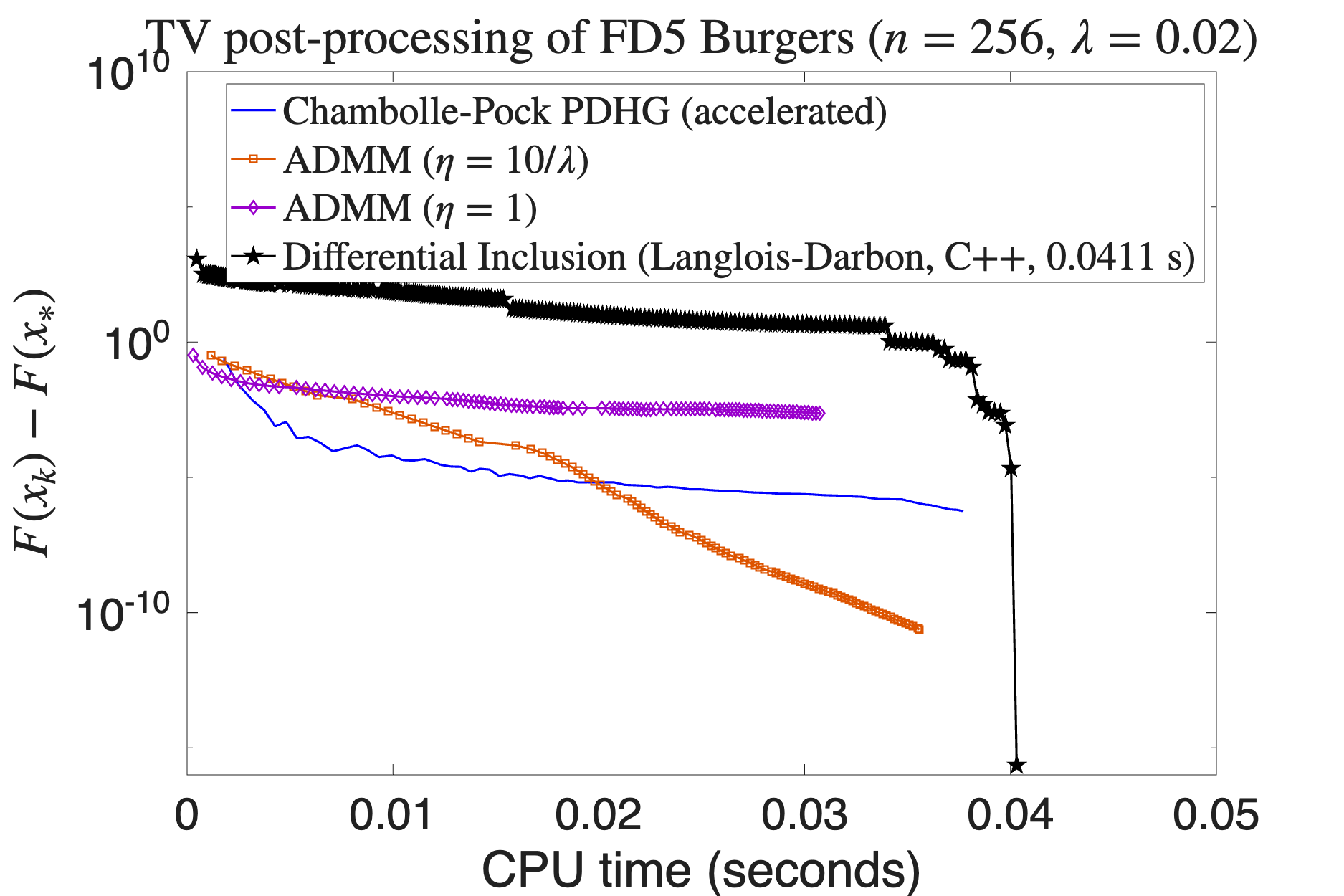}
    \caption{Objective error vs.\ CPU time.}
    \label{fig:burgers-shock-conv-time}
  \end{subfigure}
  \caption{Convergence comparison for TV post-processing of the
  FD5 Burgers solution ($n = 256$, $\lambda = 0.02$).
  The DI algorithm converges in 244~steps,
  comparable to the hardest spectral case (496~steps at $N=1000$).}
  \label{fig:burgers-shock-conv}
\end{figure}

Figure~\ref{fig:burgers-shock-conv} compares the convergence of the
three TV solvers on the post-processing task with $\lambda = 0.02$.
The DI algorithm converges in
244~breakpoint steps.
The step count is comparable to the hardest spectral denoising test
(496~steps at $N=1000$) and reflects the richer structure of the FD5 Burgers
solution: the oscillatory signal near the shock produces many
nonzero differences, so the equicorrelation set grows to a large
fraction of the $n - 1 = 255$ possible indices.
ADMM with tuned $\eta = 10/\lambda$ converges to ${\sim}\,10^{-11}$
within 500~iterations. PDHG reaches ${\sim}\,10^{-7}$.
Untuned ADMM stalls near $10^{-3}$.

\FloatBarrier
\subsection{Tests on compressible Euler equations}
\label{sec:fd-euler}

Consider the 1D compressible Euler equations:
\begin{equation}\label{eq:euler1d}
  \frac{\partial}{\partial t}
  \begin{pmatrix} \rho \\ \rho u \\ E \end{pmatrix}
  + \frac{\partial}{\partial x}
  \begin{pmatrix} \rho u \\ \rho u^2 + p \\ (E+p)\,u \end{pmatrix}
  = 0,
\end{equation}
where $\rho$ is density, $u$ is velocity, $p$ is pressure, and
$E = p/(\gamma-1) + \tfrac{1}{2}\rho u^2$ is the total energy with
$\gamma = 1.4$.

We discretize~\eqref{eq:euler1d} with the same fifth-order linear
upwind FD5 scheme from Section~\ref{sec:pde-burgers}, applied
componentwise with Lax--Friedrichs flux splitting, integrated
in time with the standard third-order SSP Runge--Kutta method.
The resulting scheme is locally conservative and fifth-order accurate in
smooth regions, but oscillatory near discontinuities.
To ensure that density
and pressure remain positive, one can augment the scheme with the
Zhang--Shu positivity-preserving (PP)
limiter for finite difference schemes~\cite{ZhangShu2012fd}, which ensures that the SSP Runge--Kutta update preserves the admissible set
$G = \{(\rho,\rho u,E)^T : \rho > 0,\; p > 0\}$.  The PP limiter
guarantees positivity but does not suppress oscillations.
We
use an SSP scheme here, rather than the RK4 of
Section~\ref{sec:pde-burgers}, because the Zhang--Shu
positivity-preserving limiter requires the time stepping to be a
convex combination of forward-Euler steps.
To suppress the oscillations, the TV solver can be used as a \emph{post-processing} step to the primitive variables $(\rho, u, p)$ after time integration is complete.
Because the limiter is applied to the primitive variables, it
preserves the discrete mean of each primitive but not the conserved
momentum $\rho u$ or total energy $E$; conservation in the sense of
Section~\ref{sec:formulation} holds only for the quantity actually
denoised. For the post-processing tests here this is acceptable, but a
conservative variant would denoise the conserved variables directly.

\begin{figure}[!htbp]
  \centering
  \includegraphics[width=\textwidth]{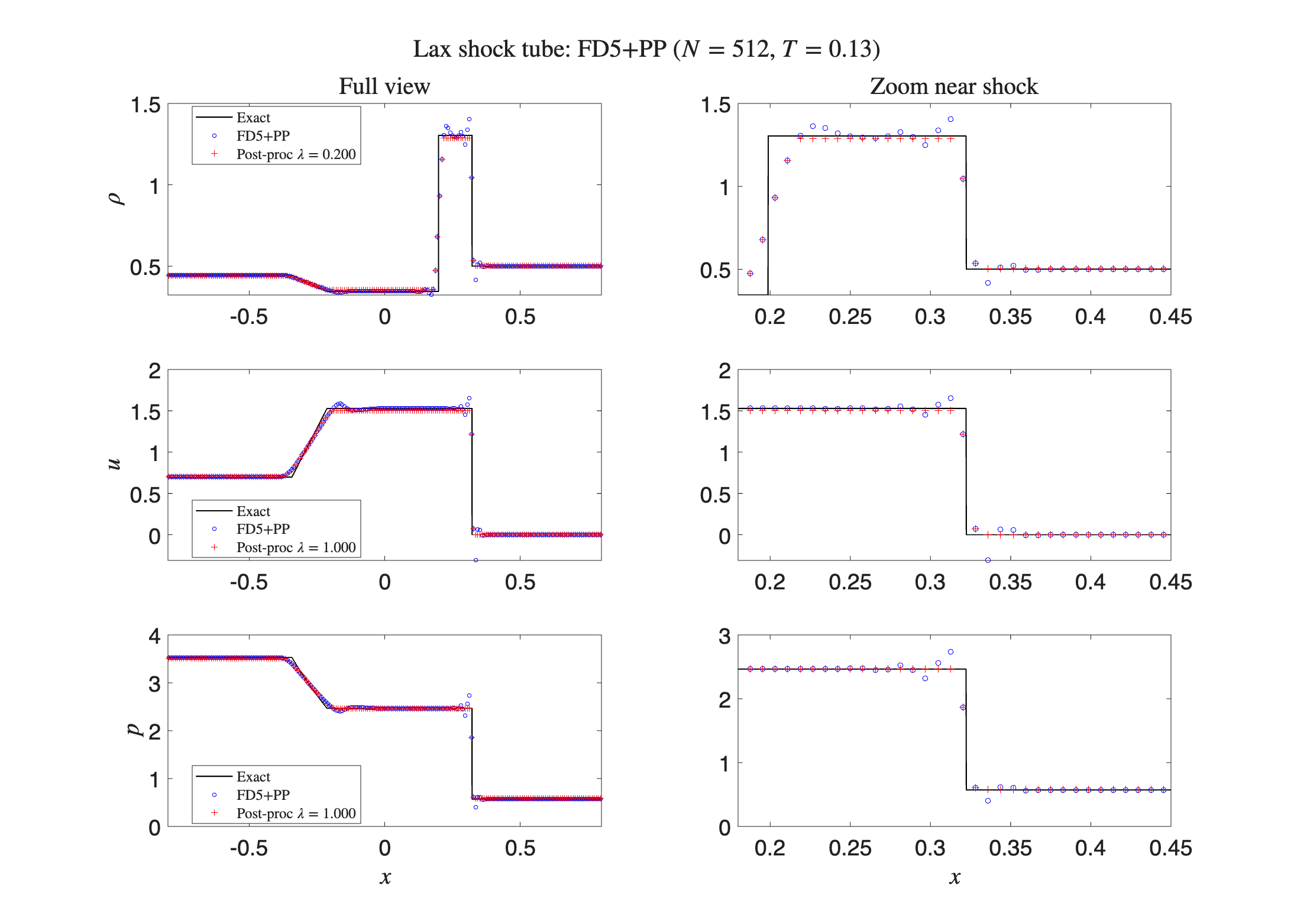}
  \caption{Lax shock tube with FD5+PP base scheme ($N = 512$, $T = 0.13$).
  Left: full view. Right: zoom near shock.
  Blue circles: FD5+PP numerical solution.
  Red plus signs: TV post-processed primitives
  ($\lambda_\rho = 0.2$, $\lambda_u = 1.0$, $\lambda_p = 1.0$).
  Black line: exact Riemann solution.}
  \label{fig:fd-lax-pp}
\end{figure}

\paragraph{Lax shock tube}
Figure~\ref{fig:fd-lax-pp} shows results for the Lax shock
tube~\cite{lax1954} ($N = 512$, $T = 0.13$).  The unprocessed
FD5+PP solution (blue circles) has $L^2(\rho) = 1.195 \times 10^{-1}$
with visible oscillations near the contact discontinuity and shock.
Because TV denoising acts on each variable independently, different
regularization strengths can target different wave features: the
density profile contains the contact discontinuity in addition to the
shock, so moderate smoothing is used to avoid flattening the contact,
whereas velocity and pressure are continuous across the contact and
tolerate stronger regularization to remove their oscillations.  With per-variable parameters $\lambda_\rho = 0.2$,
$\lambda_u = 1.0$, $\lambda_p = 1.0$ (red plus signs), TV post-processing
removes oscillations and reduces the density $L^2$ error by 4.4\%
($ 1.195 \times 10^{-1}$ to $1.143 \times 10^{-1}$).

\FloatBarrier
\subsection{Rigid body rotation: LeVeque three-body test}
\label{sec:pde-rotation}

We next solve the rigid body rotation problem
\[
  u_t - y\,u_x + x\,u_y = 0,
  \qquad (x,y) \in [-\pi,\pi)^2, \quad T = 2\pi,
\]
following Example~4.7 of Liu, Cheng, and Shu~\cite{liu2017sweep}.
The velocity field $\bv = (-y, x)$ is divergence-free and produces
solid-body rotation about the origin with period~$2\pi$.  The initial
condition is the LeVeque three-body test~\cite{leveque1996}: a cosine
bell centered at $(-\pi/2, 0)$, a cone centered at $(0, -\pi/2)$,
and a slotted disk centered at $(0, \pi/2)$, all with radius
$r_0 = 0.3\pi$.  The exact solution at $T = 2\pi$ equals the initial
condition.

We discretize on a $128 \times 128$ grid using the Fourier
pseudospectral method with Orszag's $\tfrac{2}{3}$ dealiasing
rule~\cite{orszag1971} and RK4 time-stepping
($\Delta t = 2/(\pi N) \approx 4.97 \times 10^{-3}$).  Nonlinear
products transfer energy to higher wavenumbers than are present in the
inputs; any content above the grid's Nyquist cutoff $|k| = N/2$ wraps
around (``aliases'') and pollutes the low-wavenumber spectrum.  The
$\tfrac{2}{3}$ rule restricts the solution to modes with $|k| < N/3$
(cutoff at $\tfrac{2}{3}$ of the Nyquist wavenumber $N/2$), zeroing
the remaining modes before and after each nonlinear step.  The product of two such inputs has spectrum in
$|k| < 2N/3$, and any aliased content folds into $|k| \ge N/3$, which
the filter discards.
We compare the methods in Table~\ref{tab:rotation}. The PDE solver is
implemented in MATLAB, and the reported TV time is from a MATLAB wrapper
calling the C++ DI solver.
\begin{table}[ht]
\centering\footnotesize\setlength{\tabcolsep}{3pt}
\begin{tabular}{lcccccccc}
  \hline
  Method & $\lambda$ & $L^1$  {error} & min & max & DI steps & {TV calls} & TV (s)
    & PDE (s) \\
  \hline
  {No limiter} & {---} & {1.72e-2} & {$-$0.202} & {1.257}
    & {---} & {---} & {---} & {3.5} \\
  TV limiter ($K{=}10$) & 0.02 & 1.94e-2 & $-$0.025 & 0.983
    & 76753 & {12} & 0.35 & 3.5 \\
  TV post-proc & 0.05 & 6.81e-3 & $-$0.012 & 1.057
    & 5414 & {1} & 0.02 & 3.5 \\
  TV post-proc & 0.10 & 9.14e-3 & 0.001 & 1.000
    & 4457 & {1} & 0.01 & 3.5 \\
  \hline
\end{tabular}
\caption{Rigid body rotation ($128 \times 128$, $T = 2\pi$).
Exact solution satisfies $u \in [0, 1]$.
DI steps: total breakpoint iterations across all 1D solves.
TV/PDE: wall-clock time (seconds) for TV operations and PDE solver.
Post-processing uses the no-limiter PDE solution.}
\label{tab:rotation}
\end{table}

\begin{figure}[!htbp]
  \centering
  \includegraphics[width=\textwidth]{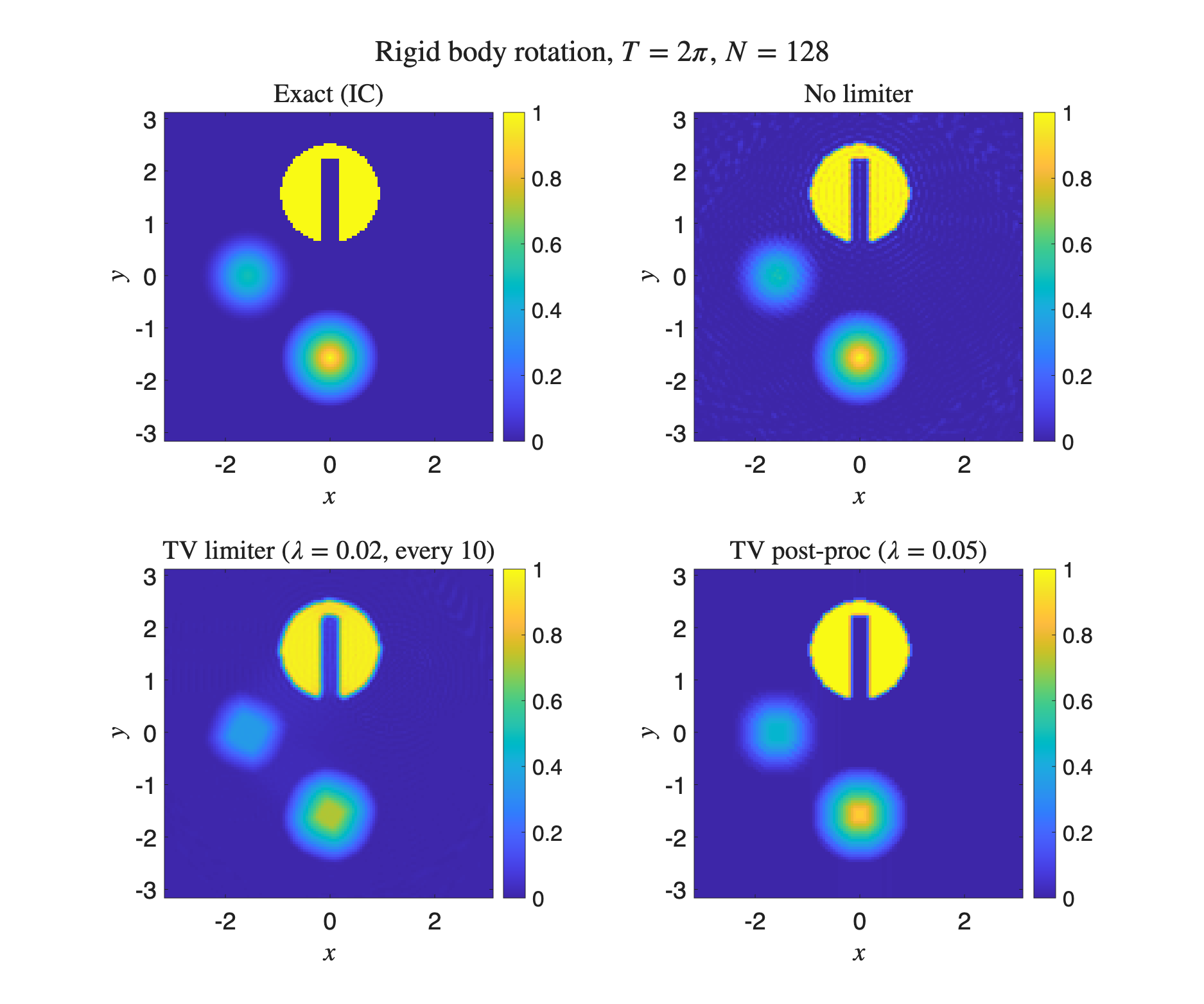}
  \caption{Rigid body rotation at $T = 2\pi$ ($128 \times 128$).
  Top left: exact (= IC).  Top right: spectral (no limiter).
  Bottom left: TV limiter ($K{=}10$, $\lambda{=}0.02$).
  Bottom right: TV post-processing ($\lambda{=}0.05$).}
  \label{fig:tv-2d-rotation}
\end{figure}

\begin{figure}[!htbp]
  \centering
  \includegraphics[width=\textwidth]{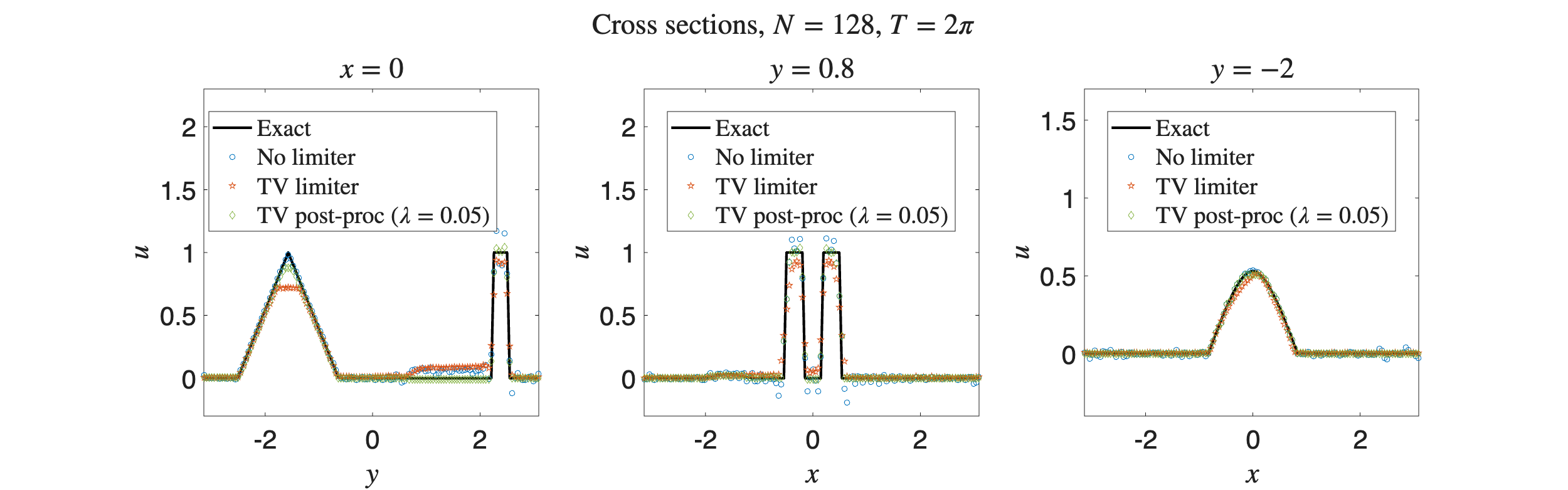}
  \caption{Cross-sectional cuts for the rotation test.
  Left: $x = 0$ (through cone and slotted disk).
  Center: $y = 0.8$ (through slotted disk).
  Right: $y = -2$ (near the cosine bell).
  Solid black line: exact solution.
  Circles (blue): no limiter. Stars (orange): TV limiter ($K{=}10$).
  Diamonds (green): TV post-processing ($\lambda{=}0.05$).
  Both TV approaches significantly reduce the Gibbs oscillations
  visible in the no-limiter solution, with post-processing providing
  the best $L^1$ accuracy.}
  \label{fig:tv-2d-rotation-cuts}
\end{figure}

Without any limiter, the spectral solution overshoots by 25.7\% and
undershoots by 20.2\%.  The per-step TV limiter reduces bound
violations to $[-0.025, 0.983]$ with only 12~TV calls over 1264~time
steps.  TV post-processing at $\lambda = 0.10$ achieves near-perfect bounds
$[0.001, 1.000]$ with the smallest undershooting.  The best $L^1$
accuracy ($6.81 \times 10^{-3}$) is obtained by post-processing with
$\lambda = 0.05$, which balances oscillation removal against over-smoothing. We note that the per-step limiter slightly increases the
$L^1$ error relative to the no-limiter solution (from
$1.72 \times 10^{-2}$ to $1.94 \times 10^{-2}$): repeated denoising
accumulates smoothing error, which is the price paid for controlling
the oscillations during the computation. When only the final solution
is needed, post-processing avoids this accumulation. The full solutions and the corresponding cross-sectional cuts are shown
in Figures~\ref{fig:tv-2d-rotation} and~\ref{fig:tv-2d-rotation-cuts}.

\FloatBarrier
\subsection{Swirling deformation flow}
\label{sec:pde-swirl}

The swirling deformation problem has a time-dependent,
spatially varying velocity field (Example~4.8 of~\cite{liu2017sweep}):
\[
  u_t - \Bigl(\cos^2\!\bigl(\tfrac{x}{2}\bigr)\sin(y)\,g(t)\,u\Bigr)_x
      + \Bigl(\sin(x)\cos^2\!\bigl(\tfrac{y}{2}\bigr)\,g(t)\,u\Bigr)_y = 0,
\]
where $g(t) = \pi\cos(\pi t/T)$ and $T = 1.5$.  The velocity field
is divergence-free. It deforms the solution into thin filaments during
$0 < t < T/2$, then reverses, so the exact solution at $T = 1.5$
equals the initial condition.  The IC is the same LeVeque three-body
test as in Section~\ref{sec:pde-rotation}.

\begin{figure}[!htbp]
  \centering
  \includegraphics[width=\textwidth]{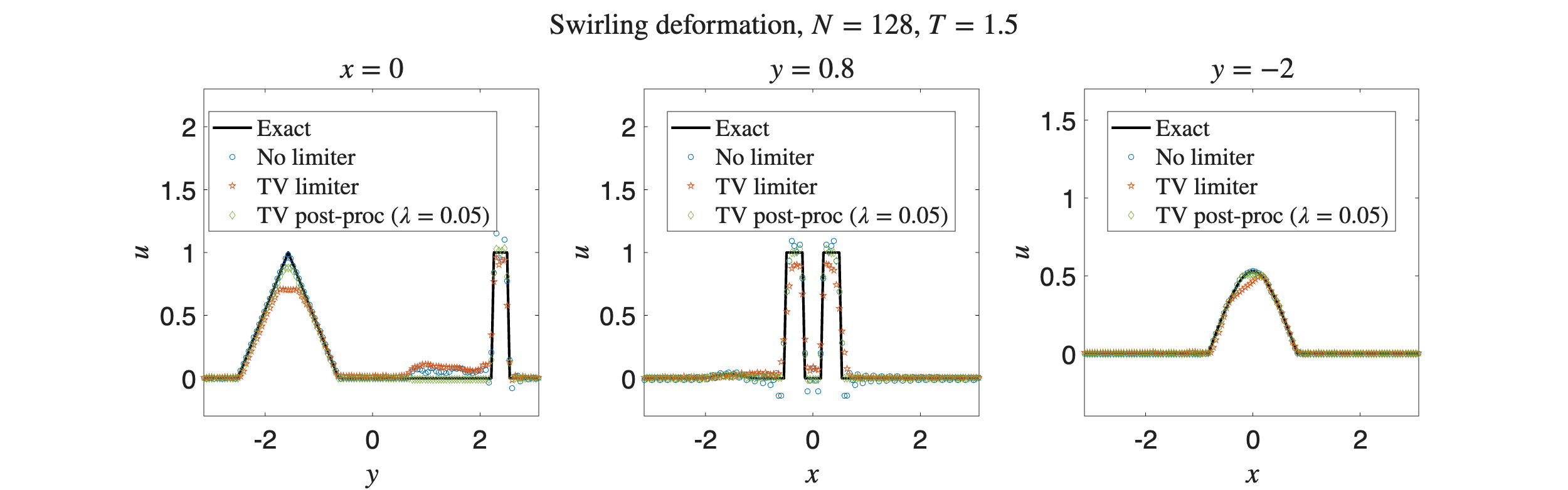}
  \caption{Cross-sectional cuts for the swirling deformation test.
  Same cuts and marker conventions as Figure~\ref{fig:tv-2d-rotation-cuts}.
  The TV limiter and post-processing effectively suppress oscillations
  despite the time-dependent, spatially varying velocity field.}
  \label{fig:tv-2d-swirl-cuts}
\end{figure}

On a $128 \times 128$ grid, the TV limiter reduces bound violations from $[-0.19, 1.25]$
to $[-0.035, 1.008]$ with 13~TV calls over 302~time steps.
Post-processing at $\lambda = 0.05$ gives the best $L^1$ accuracy
($6.20 \times 10^{-3}$), while $\lambda = 0.10$ achieves near-perfect
bounds $[0.001, 0.994]$. As in Section~\ref{sec:pde-rotation}, the
per-step limiter trades some $L^1$ accuracy for improved bounds. The corresponding cross-sectional cuts through
the three bodies are shown in Figure~\ref{fig:tv-2d-swirl-cuts}.

\begin{table}[ht]
\centering
\begin{tabular}{lcccc}
  \hline
  Method & $L^1$ error & min & max & TV calls \\
  \hline
  No limiter & $1.03 \times 10^{-2}$ & $-0.188$ & $1.249$ & --- \\
  TV limiter ($K{=}10$, $\lambda{=}0.02$)
    & $2.21 \times 10^{-2}$ & $-0.035$ & $1.008$ & 13 \\
  TV post-proc ($\lambda{=}0.05$)
    & $6.20 \times 10^{-3}$ & $-0.010$ & $1.049$ & 1 \\
  TV post-proc ($\lambda{=}0.10$)
    & $8.65 \times 10^{-3}$ & $0.001$ & $0.994$ & 1 \\
  \hline
\end{tabular}
\caption{Swirling deformation flow ($128 \times 128$, $T = 1.5$).
Exact solution satisfies $u \in [0, 1]$.}
\label{tab:swirl}
\end{table}

\FloatBarrier
\subsection{Incompressible Euler: vortex patches}
\label{sec:pde-euler}

We finally consider the incompressible 2D Euler equations in
vorticity form (Example~4.10 of~\cite{liu2017sweep}):
\[
  \omega_t + u\,\omega_x + v\,\omega_y = 0,
  \qquad \Delta\psi = \omega,
  \qquad (u,v) = (-\psi_y, \psi_x),
\]
on $[0,2\pi)^2$ with periodic boundary conditions.  The initial
vorticity consists of two rectangular patches:
$\omega_0 = -1$ on $[\pi/2, 3\pi/2] \times [\pi/4, 3\pi/4]$,
$\omega_0 = 1$ on $[\pi/2, 3\pi/2] \times [5\pi/4, 7\pi/4]$,
and $\omega_0 = 0$ elsewhere.  The maximum principle guarantees
$\omega \in [-1, 1]$ at all times.

We discretize on a $128 \times 128$ grid using the Fourier
pseudospectral method ($\tfrac{2}{3}$ dealiasing, RK4) and integrate
to $T = 5$.  The velocity is recovered at each stage by solving the
Poisson equation $\Delta\psi = \omega$ in Fourier space.
Table~\ref{tab:euler} summarizes the results. The PDE solver is
implemented in MATLAB, and the reported TV time is from a MATLAB wrapper
calling the C++ DI solver.

\begin{table}[ht]
\centering\footnotesize\setlength{\tabcolsep}{4pt}
\begin{tabular}{lccccccc}
  \hline
  Method & $\lambda$ & min & max & DI steps & {TV calls} & TV (s)
    & PDE (s) \\
  \hline
  {No limiter} & {---} & {$-$1.299} & {1.299}
    & {---} & {---} & {---} & {1.9} \\
  TV limiter ($K{=}10$) & 0.02 & $-$0.945 & 0.945
    & 265119 & {41} & 0.95 & 1.7 \\
  TV post-proc & 0.10 & $-$1.053 & 1.053
    & 3631 & {1} & 0.01 & 1.9 \\
  \hline
\end{tabular}
\caption{Incompressible Euler ($128 \times 128$, $T = 5$).
Maximum principle: $\omega \in [-1, 1]$.
DI steps: total breakpoint iterations across all 1D solves.
TV/PDE: wall-clock time (seconds) for TV operations and PDE solver.
Post-processing uses the no-limiter PDE solution.}
\label{tab:euler}
\end{table}

\begin{figure}[!htbp]
  \centering
  \includegraphics[width=\textwidth]{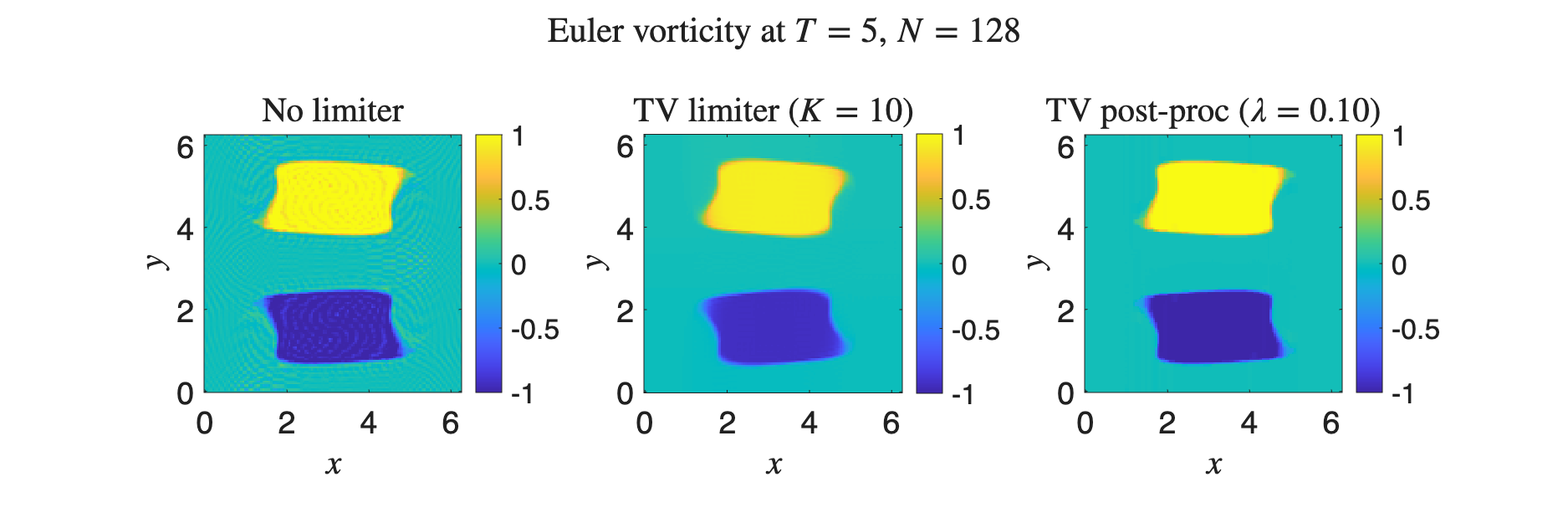}
  \caption{Vorticity fields at $T = 5$ ($128 \times 128$).
  Left: no limiter.  Center: TV limiter ($K{=}10$, $\lambda{=}0.02$).
  Right: TV post-processing ($\lambda{=}0.10$).}
  \label{fig:tv-2d-euler-vort}
\end{figure}

\begin{figure}[!htbp]
  \centering
  \includegraphics[width=\textwidth]{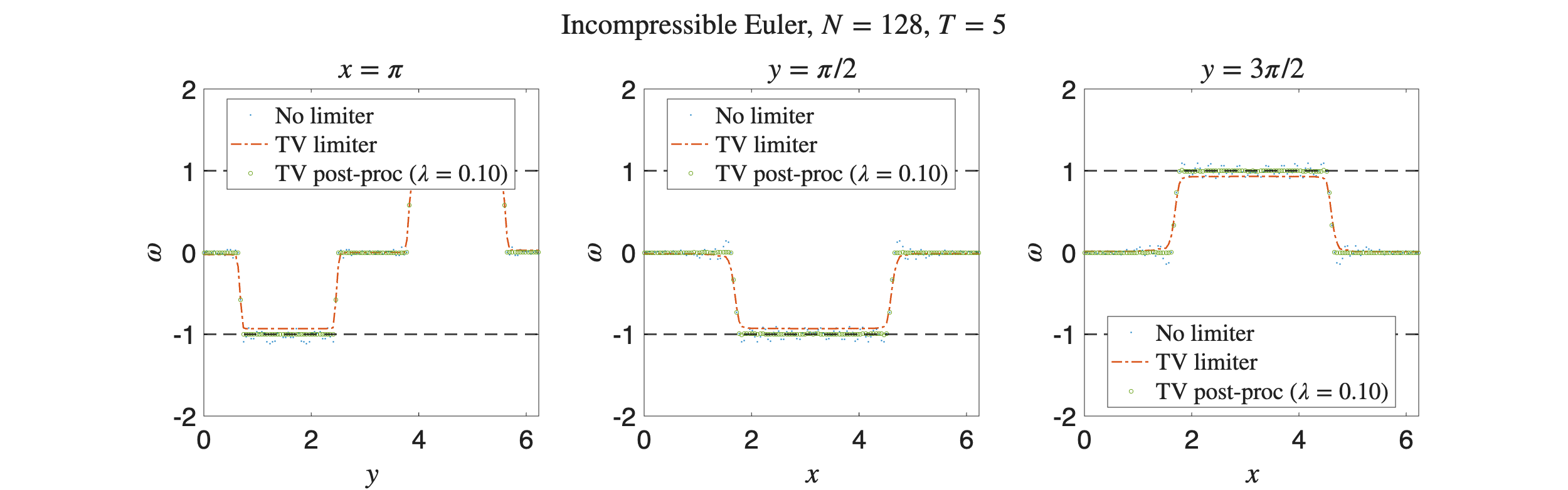}
  \caption{Cross-sectional cuts for incompressible Euler at $T = 5$.
  Left: $x = \pi$.  Center: $y = \pi/2$.
  Right: $y = 3\pi/2$.
  Dots (blue): no limiter. Dash-dot (orange): TV limiter ($K{=}10$).
  Circles (green): TV post-processing ($\lambda{=}0.10$).
  Dashed lines mark the maximum principle bounds $\omega = \pm 1$.}
  \label{fig:tv-2d-euler-cuts}
\end{figure}

Without any limiter, the spectral solution violates the maximum
principle by 30\%.  The per-step TV limiter reduces the range to
$[-0.945, 0.945]$, well within the physical bounds, with 41~TV calls
over 435~time steps.  TV post-processing at $\lambda = 0.10$ reduces the
violation to about 5\%, but cannot fully recover the bounds from the
accumulated nonlinear oscillations.
The vorticity fields and the corresponding cross-sectional cuts are
shown in Figures~\ref{fig:tv-2d-euler-vort}
and~\ref{fig:tv-2d-euler-cuts}.

\FloatBarrier
\subsection{2D Riemann problem with FD5 scheme}
\label{sec:pde-riemann2d}

We also solve a genuine two-dimensional Riemann problem
for the compressible Euler equations using the fifth-order linear
finite difference scheme (FD5) from Section~\ref{sec:pde-burgers},
extended to two dimensions with Lax--Friedrichs splitting applied
dimension by dimension.  We use Configuration~6 of Lax and
Liu~\cite{laxliu1998}, in which the initial data consist of four
constant states meeting at the center of the domain $[0,1]^2$:
\[
  (\rho, u, v, p) =
  \begin{cases}
    (2, 0.75, 0.5, 1) & x < 0.5,\; y \geq 0.5, \\
    (1, 0.75, {-}0.5, 1) & x \geq 0.5,\; y \geq 0.5, \\
    (1, {-}0.75, 0.5, 1) & x < 0.5,\; y < 0.5, \\
    (3, {-}0.75, {-}0.5, 1) & x \geq 0.5,\; y < 0.5.
  \end{cases}
\] 
We compute a reference solution with the WENO5
scheme~\cite{JiangShu1996} on an $800 \times 800$ grid with SSP-RK3
time stepping.

\begin{figure}[!htbp]
  \centering
  \includegraphics[width=\textwidth]{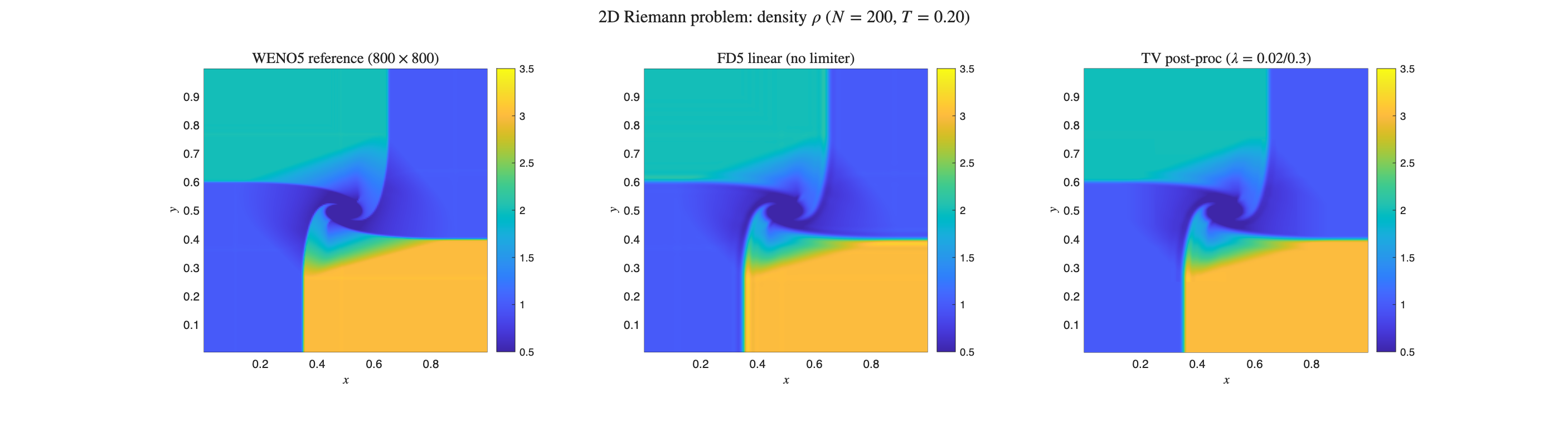}
  \caption{Density $\rho$ for the 2D Riemann problem (Config.~6,
  $200 \times 200$, $T = 0.2$).
  Left: WENO5 reference ($800 \times 800$).
  Center: FD5 linear (no limiter).
  Right: TV post-processing ($\lambda = 0.02$ center, $0.3$ corners).}
  \label{fig:riemann2d-density}
\end{figure}

\begin{figure}[!htbp]
  \centering
  \includegraphics[width=\textwidth]{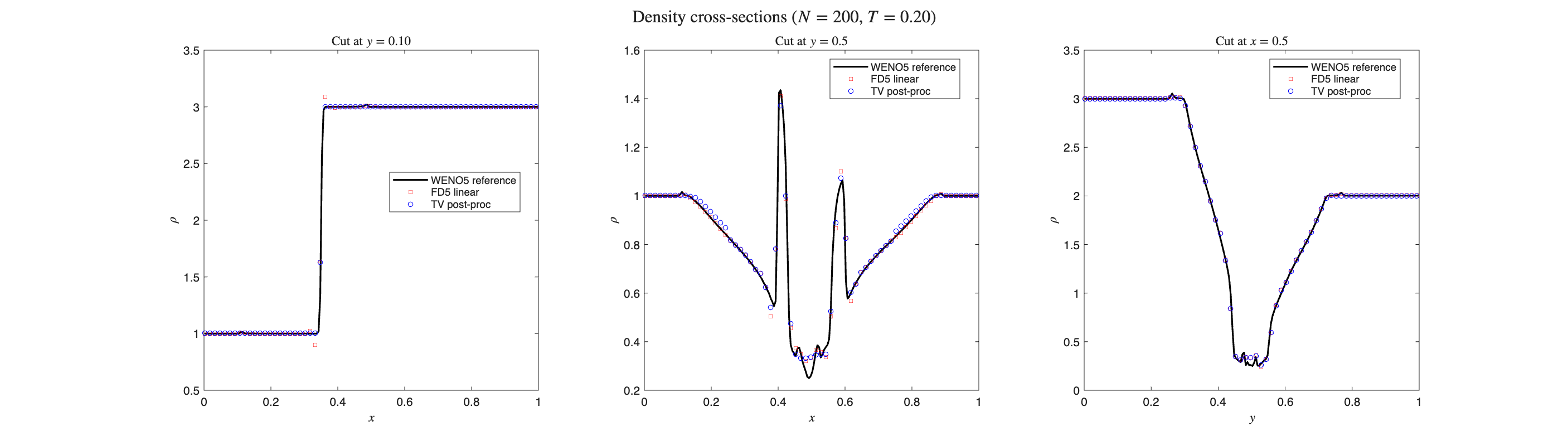}
  \caption{Density cross-sections for the 2D Riemann problem.
  Left: cut at $y = 0.10$.  Center: cut at $y = 0.5$.
  Right: cut at $x = 0.5$.
  Black solid: WENO5 reference ($800 \times 800$).
  Red squares: FD5 linear. Blue circles: TV post-processing ($\lambda = 0.02/0.3$).}
  \label{fig:riemann2d-cuts}
\end{figure}

Both the FD5 and WENO5 schemes use zero-order extrapolation
(outflow) boundary conditions, following~\cite{laxliu1998}.
We discretize on a $200 \times 200$ grid with RK4 time stepping
($\mathrm{CFL} = 0.5$) and integrate to $T = 0.2$.
Figure~\ref{fig:riemann2d-density} compares the density fields. The
WENO5 reference (left) resolves sharp shock and contact structures,
while the FD5 linear scheme (center) captures the wave structure but
produces oscillations near all discontinuities.
For TV post-processing (right), we use a per-slice regularization
parameter.  Since the 2D TV solver operates dimension by dimension
(columns then rows), each 1D slice receives its own~$\lambda$.
Slices passing through the center of the domain
($x \in [0.25, 0.75]$ for columns, $y \in [0.25, 0.75]$ for rows)
use $\lambda = 0.02$, preserving the complex but not strongly
oscillatory wave interactions. Slices in the corner regions use
$\lambda = 0.3$, providing stronger smoothing where the solution
is nearly constant.
TV post-processing reduces the oscillations near these
discontinuities in the cross-sectional cuts through $y = 0.10$,
$y = 0.5$, and $x = 0.5$. See Figure~\ref{fig:riemann2d-cuts}.

The PDE solver requires 358~time steps (12.3\,s).  TV post-processing
of the density field takes 19107~DI steps (0.46\,s), at a fraction of
the PDE solver cost.

\FloatBarrier
\section{Concluding remarks}
\label{sec:conclusions}

In this work, we applied the DI algorithm to the one-dimensional  TV denoising
problem and used it as a non-oscillatory limiter for high-order methods
for conservation laws. The key idea is that the difference operator $D \in \R^{(n-1)\times n}$ is surjective, so the change of variables
$\bz = D\tilde\bx$ reduces the TV problem
to an equivalent LASSO problem, which the DI algorithm can solve exactly in finitely many steps. As a result, we obtain a 1D solver that computes the exact TV
minimizer and preserves the total mass
of the input. The only free parameter is $\lambda$, which directly
controls the regularization strength. In two dimensions, we use a
dimension-by-dimension splitting.  We have tested the TV limiter on Fourier pseudospectral methods and a fifth-order finite difference scheme solving 1D and 2D problems. 

\appendix
\section{Comparison of analysis operators}
\label{app:operators}

The denoising framework generalizes to any analysis operator~$D$:
\begin{equation}\label{eq:gen-denoise}
  \min_{\bx}\;\norm{D\bx}_1 + \frac{1}{2\lambda}\norm{\bx - \bb}_{2}^{2}.
\end{equation}
The TV problem studied in this paper corresponds to choosing $D$ as the
forward difference operator.  As observed by Cai et
al.~\cite{cai2012image}, TV denoising is in fact the simplest member of
a family of B-spline framelet models.  The Haar framelet (Example~2.1
of~\cite{cai2012image}), constructed from the piecewise constant
B-spline $B_1$, has a single detail mask
$\ba_1 = \tfrac{1}{2}[1,-1]$, which is the forward difference operator
up to a factor of~$1/2$.  Consequently, the Haar framelet
analysis-based approach and the Rudin--Osher--Fatemi (ROF)
model~\cite{rudin1992} coincide under forward difference
discretization~\cite[eq.~(2.3)]{cai2012image}.  Higher-order
B-spline framelets add detail bands that approximate higher-order
difference operators, yielding natural generalizations of TV.

The piecewise linear framelet (Example~2.2 of~\cite{cai2012image}),
constructed from the B-spline $B_2$, has two detail bands: $\ba_1$
approximates a first-order difference and $\ba_2$ a second-order
difference ($2n\times n$ tight frame).  The piecewise cubic framelet
(Example~2.3 of~\cite{cai2012image}), constructed from $B_4$, has four detail bands
approximating first- through fourth-order differences ($4n\times n$
tight frame).  In all cases, the analysis operator $D$ is formed by
stacking the detail-band circulant matrices (excluding the low-pass
band), and its null space is one-dimensional (the low-pass mode).

In addition, we compare against multi-level orthogonal discrete wavelet
transforms (DWT).  For an $n\times n$ orthogonal matrix~$W$ (e.g.\
Daubechies, Coiflet, or Symlet), the substitution $\bz = W\bx$
decouples the problem into componentwise soft thresholding:
$\bz^* = S_\lambda(W\bb)$, $\bx^* = W^T\bz^*$.
Equivalently, the DI algorithm applied to the LASSO~\eqref{eq:lasso-final}
with measurement matrix $M = I$ and
data $W\bb$ in place of $\tilde\bb$
recovers the same answer in a small number of
DI steps (2--30).

The LASSO reduction of \S\ref{sec:reduction} is exact only when $D$
is surjective (the forward-difference operator) or square and invertible
(the orthogonal DWTs). For the redundant B-spline framelets ($D$ tall,
$m > n$), setting $M = D^\dagger$ and minimizing
$\norm{\bz}_1 + \frac{1}{2\lambda}\norm{M\bz - \tilde{\bb}}_2^2$ over
$\bz \in \R^m$ solves the \emph{synthesis} formulation; the recovered
$\bx^* = M\bz^*$ minimizes a relaxation of the analysis
problem~\eqref{eq:gen-denoise}, not~\eqref{eq:gen-denoise} itself, since
the minimizer need not satisfy $\bz^* = D\bx^*$. We report the synthesis
solutions below; they coincide with the analysis minimizer only in the
surjective and orthogonal cases.  We test the following operators on $n = 256$
grid points using the same two Gibbs oscillation test
signals from~\S\ref{sec:tv-spectral} and~\S\ref{sec:tv-spectral2}, with $\lambda = 0.5$ for the
piecewise constant signal and $\lambda = 0.1$ for the piecewise smooth
signal:
\begin{enumerate}[nosep]
\item \emph{Forward difference} (TV / Haar framelet): $(n{-}1)\times n$,
  single detail band (first-order difference).
\item \emph{Multi-level orthogonal DWT}: $n\times n$,
  db1 (Haar) and db4 (Daubechies-4).
\item \emph{Piecewise linear framelet}: $2n\times n$,
  two detail bands (first- and second-order differences).
\item \emph{Piecewise cubic framelet}: $4n\times n$,
  four detail bands (first- through fourth-order differences).
\end{enumerate}
In every case the DI algorithm solves its
LASSO~\eqref{eq:lasso-final} exactly after building $M = D^\dagger$
offline (or $M = I$ for the orthogonal DWT case); for the redundant
framelets this LASSO is the synthesis relaxation rather than the
analysis problem~\eqref{eq:gen-denoise}, as noted above.

\begin{figure}[!htbp]
  \centering
  \includegraphics[width=\textwidth]{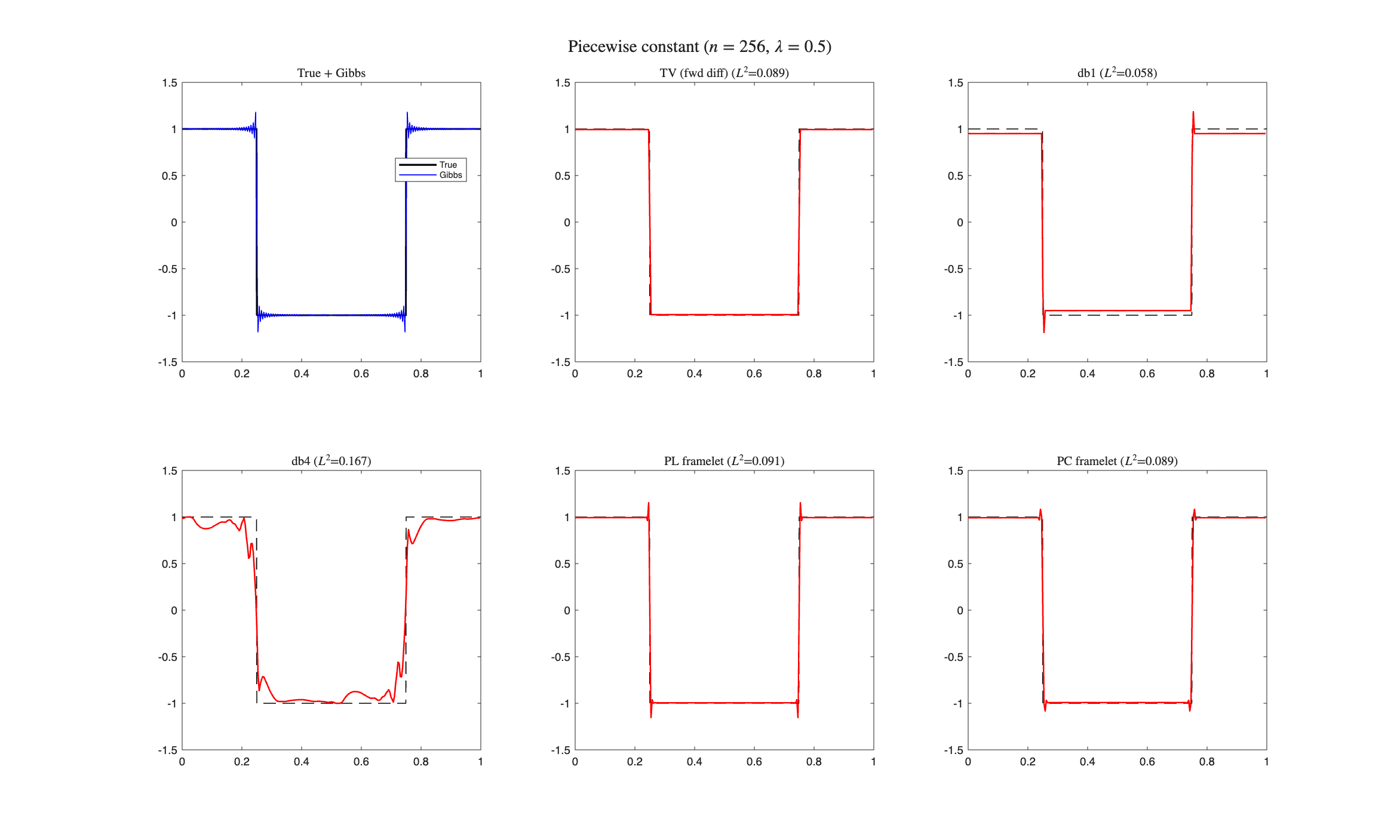}
  \caption{Denoised piecewise constant signal ($n = 256$, $\lambda = 0.5$).
    TV and  non-TV operators are shown.}
  \label{fig:op-pwconst}
\end{figure}

\begin{figure}[!htbp]
  \centering
  \includegraphics[width=\textwidth]{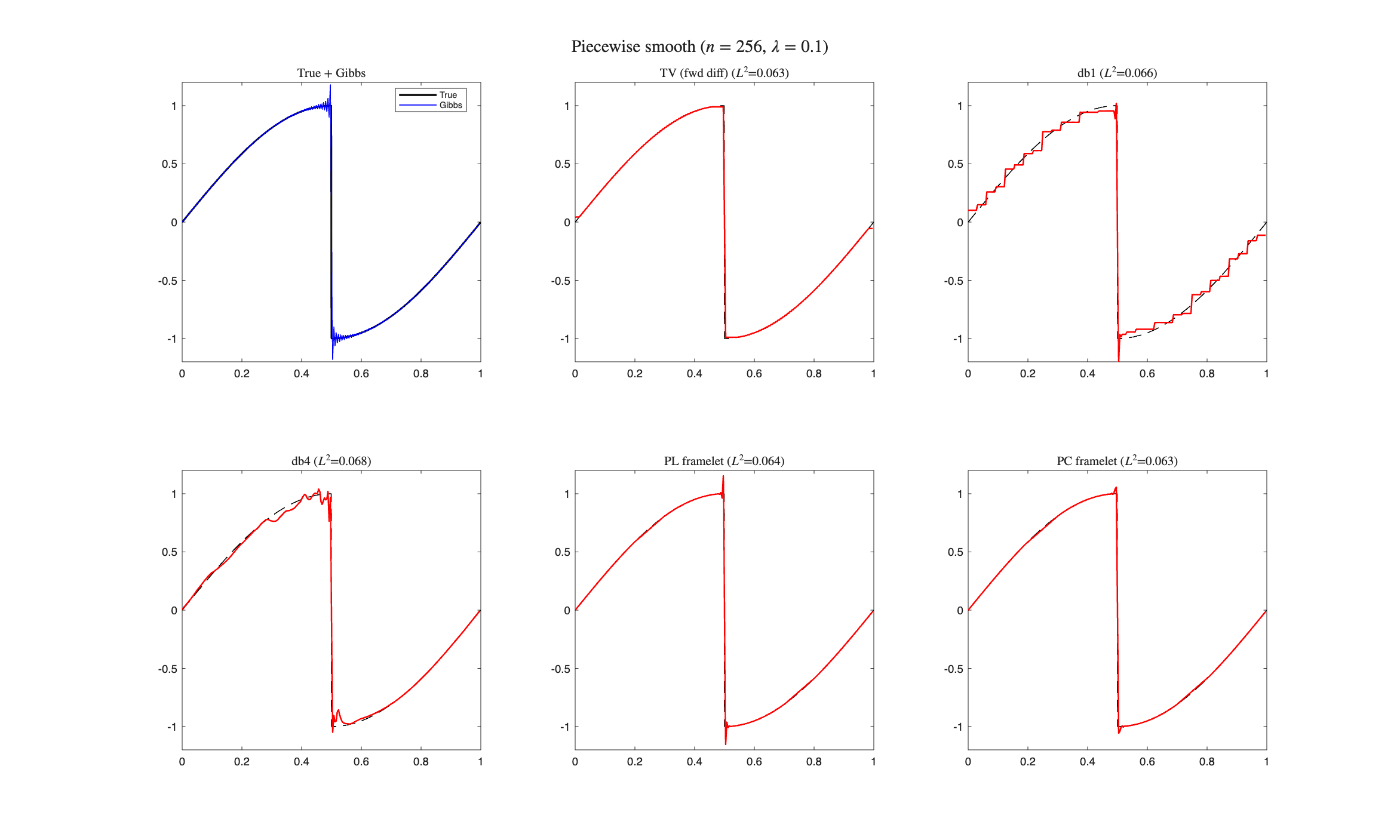}
  \caption{Denoised piecewise smooth signal ($n = 256$, $\lambda = 0.1$).
    TV and  non-TV operators are shown.}
  \label{fig:op-pwsmooth}
\end{figure}

Table~\ref{tab:operators} reports the number of DI steps and $L^2$
error for each operator.
Figures~\ref{fig:op-pwconst}~and~\ref{fig:op-pwsmooth} show the
denoised signals for all operators.

\begin{table}[!htbp]
  \centering
  \caption{Comparison of analysis operators ($n = 256$). The DI algorithm
    solves the associated LASSO exactly (the synthesis form for the
    redundant framelets).
    Entries are piecewise constant ($\lambda = 0.5$) / piecewise smooth ($\lambda = 0.1$).}
  \label{tab:operators}
  \smallskip
  \begin{tabular}{llcc}
    \hline
    Operator & Type & DI steps &
      $\norm{\bx^*\!-\!\bx_{\rm true}}_2$ \\
    \hline
    Forward diff.\ (TV)
      & Haar framelet & 3 / 139 & 0.089 / 0.063 \\
    db1 (multi-level Haar)
      & $n{\times}n$ orth.\ DWT & 2 / 30 & 0.058 / 0.066 \\
    db4
      & $n{\times}n$ orth.\ DWT & 11 / 24 & 0.167 / 0.068 \\
    PL framelet
      & $2n{\times}n$ tight frame & 3 / 148 & 0.091 / 0.064 \\
    PC framelet
      & $4n{\times}n$ tight frame & 1 / 773 & 0.089 / 0.063 \\
    \hline
  \end{tabular}
\end{table}

For the piecewise constant signal ($\lambda = 0.5$), the multi-level Haar DWT
(db1) achieves the lowest $L^2$ error ($0.058$), ahead of TV ($0.089$)
and all other operators. However, the smaller $L^2$ error does not fully
reflect the solution quality. The db1 solution develops staircase
artifacts near the jumps (see Figure~\ref{fig:op-pwconst}) whereas TV
preserves sharper transitions. The $L^2$ metric is sensitive to the
single grid point where TV passes through zero, but it does not capture
this difference well. The longer-filter db4 wavelet performs worse
($L^2 = 0.17$) because its wider support smears the jumps.

For the piecewise smooth signal ($\lambda = 0.1$), all operators give
comparable $L^2$ errors ($0.063$--$0.068$), with TV slightly smaller.
Figure~\ref{fig:op-pwsmooth} shows that TV and the B-spline framelets
preserve the smooth arcs and the jump more accurately, while the
multi-level wavelets introduce oscillations or excessive rounding near
the discontinuity.

These numerical findings are consistent with the theory of the Gibbs
phenomenon in wavelet expansions. Kelly~\cite{kelly1996gibbs} proved that the Haar wavelet does \emph{not}
exhibit the Gibbs phenomenon, via an
if-and-only-if criterion on its reproducing kernel
$K_0(x,y) = \sum_k \varphi(x-k)\varphi(y-k)$: a Gibbs effect occurs near
a jump if and only if $\int_0^\infty K_0(a,u)\,du > 1$ for some
$a > 0$. For the Haar system this integral equals~$1$, and the kernel
is nonnegative---analogous to the Fej\'er kernel---so the expansion
converges without overshoot at jump discontinuities.
Shim and Volkmer~\cite{shim1996gibbs} showed, conversely, that if the scaling
function~$\varphi$ is continuous, is differentiable at some dyadic
rational~$d$ with $\varphi'(d) \neq 0$, and decays as
$|\varphi(x)| \leq C(1+|x|)^{-\beta}$ with $\beta > 3$, then the wavelet expansion exhibits a
Gibbs phenomenon on at least one side of a jump.  This covers the $C^1$ compactly supported Daubechies
wavelets ($\mathrm{db}N$, $N \geq 3$) and the $C^1$ Symlets and Coiflets
(see~\cite{raeen2008gibbs} for a detailed exposition).  Since the forward difference operator in TV denoising is
the detail mask of the Haar framelet (up to a
factor of~$1/2$), TV is expected to inherit
the Gibbs-free behavior of the Haar system; this is consistent with its
superior performance for post-processing discontinuous solutions.
We note, however, that  this connection is heuristic.

The DI algorithm applies unchanged to any analysis
operator~$D$ after forming $M = D^\dagger$ offline.
For the $n\times n$ orthogonal wavelets, the problem reduces to
componentwise soft thresholding, so there is no coupling between
coefficients. This limits their effectiveness compared with the
rectangular operators ($m \neq n$), where the DI algorithm
exploits dependencies in the analysis domain.
Among the operators tested here, TV gives the most accurate overall
reconstruction of the Gibbs-oscillatory signals. It preserves sharp
jumps and retains the smooth parts of the signal well.

\section*{Funding}
JD was supported in part by the Center for Information Geometric Mechanics and Optimization (CIGMO), a PSAAP-IV Focused Investigatory Center funded by the U.S.\ Department of Energy, National Nuclear Security Administration under Award Number DE-NA0004261.
R. Lai is partially supported by NSF grant DMS-2401297.
C.-W. Shu is partially supported by NSF grant DMS-2309249.
X. Zhang is partially supported by NSF grant DMS-2208518.

\section*{Declaration of competing interest}
The authors declare that they have no known competing financial
interests or personal relationships that could have appeared to
influence the work reported in this paper.

\section*{Data availability}
Data and code are available from the authors upon request.

\section*{CRediT authorship contribution statement}
\textbf{Gabriel P. Langlois:} Methodology, Software, Investigation, Writing -- review \& editing.
\textbf{J\'er\^ome Darbon:} Conceptualization, Methodology, Software, Writing -- review \& editing, Funding acquisition.
\textbf{Rongjie Lai:} Conceptualization, Methodology, Funding acquisition.
\textbf{Chi-Wang Shu:} Conceptualization, Methodology, Funding acquisition.
\textbf{Xiangxiong Zhang:} Conceptualization, Methodology, Software, Investigation, Writing -- original draft, Writing -- review \& editing, Supervision, Funding acquisition.

\section*{Declaration of generative AI and AI-assisted technologies in the manuscript preparation process}
During the preparation of this work, the authors used Claude Code (Anthropic) in order to assist with mathematical discussions, numerical implementation, and drafting of the manuscript. After using this tool, the authors reviewed and edited the content as needed and take full responsibility for the content of the published article.

\bibliographystyle{elsarticle-num}
\bibliography{refs}

\end{document}